\documentclass[11pt]{article}

\usepackage[utf8]{inputenc}

\usepackage{color}
\usepackage{latexsym}
\usepackage{dsfont}
\usepackage{amssymb}
\usepackage{comment}
\usepackage{graphicx}
\usepackage{amsmath,amsfonts,amssymb,theorem,euscript,array,enumerate,amsfonts,mathrsfs}
\usepackage{hyperref}
\usepackage{appendix}
\usepackage[T1]{fontenc}
\usepackage{babel}
\numberwithin{equation}{section}
\usepackage{bbm}
\usepackage{color}
\usepackage{stmaryrd} 

\usepackage{float}
\usepackage{caption}
\usepackage{subcaption}

\usepackage[normalem]{ulem}

\usepackage{comment}

\usepackage{tikz}
\usetikzlibrary{fit,matrix,chains,positioning,decorations.pathreplacing,arrows}

\usepackage{mathtools}
\mathtoolsset{showonlyrefs}

\def \trans{^{\scriptscriptstyle{\intercal}}}

\newcommand{\independent}{\protect\mathpalette{\protect\independenT}{\perp}}
\def\independenT#1#2{\mathrel{\rlap{$#1#2$}\mkern2mu{#1#2}}}

\def \b1{\bf{1}}

\def \N{\mathbb{N}}
\def \R{\mathbb{R}}

\def \E{\mathbb{E}}

\def \P{\mathbb{P}}

\def \lb{\llbracket}
\def  \rb{\rrbracket}

\def \bp{\boldsymbol{p}} 
\def \bP{\boldsymbol{P}} 
\def \bq{\boldsymbol{q}}

\def\boPhi{{\boldsymbol \Phi}}

\def \bC{\boldsymbol{C}}

\def\esssup_#1{\underset{#1}{\mathrm{ess\,sup\, }}}

\def\argmin_#1{\underset{#1}{\mathrm{argmin\, }}}
\def\argmax_#1{\underset{#1}{\mathrm{argmax\, }}}

\def \Cc{{\cal C}}

\def \Fc{{\cal F}}

\def \Ic{{\cal I}}
\def \Kc{{\cal K}}

\def \Oc{{\cal O}}
\def \Sc{{\cal S}}

\def \eps{\varepsilon}

\def \ep{\hbox{ }\hfill$\Box$}

\def\bX{{\bf X}}
\def\bY{{\bf Y}}
\def\boX{{\boldsymbol X}}
\def\boY{{\boldsymbol Y}}
\def\bolx{{\boldsymbol x}}
\def\by{{\boldsymbol y}}
\def\bolc{{\boldsymbol c}}
\def\boC{{\boldsymbol C}}
\def\boCc{{\boldsymbol \Cc}}
\def\bx{{\bf x}}
\def\bz{{\boldsymbol z}}

\def\bc{{\bf c}}

\def \mrs{\mathrm{s}}

\def \mrS{\mathrm{S}}

\def \trans{^{\scriptscriptstyle{\intercal}}}
 
\def\beqs{\begin{eqnarray*}}
\def\enqs{\end{eqnarray*}}
\def\beq{\begin{eqnarray}}
\def\enq{\end{eqnarray}}

\newcommand{\bl}[1]{\textcolor{blue}{#1}}

\newtheorem{Theorem}{Theorem}[section]
\newtheorem{Definition}{Definition}[section]
\newtheorem{Proposition}{Proposition}[section]
\newtheorem{Assumption}{Assumption}[section]

\newtheorem{Remark}{Remark}[section]
\newtheorem{Example}{Example}[section]
\numberwithin{equation}{section}

\newcommand{\dc}[1]{{\color{magenta}(DC: #1)}}

\numberwithin{equation}{section} 

\usepackage[backend=bibtex, maxbibnames=5, style=numeric,maxbibnames=99]{biblatex} 

\title{Opinion dynamics in communities with major influencers 
and implicit social influence via mean-field approximation}

\author{
Delia COCULESCU\footnote{University of Zurich, Institute for Mathematics and Department of Banking and Finance, delia.coculescu@math.uzh.ch}
\qquad\quad
M\'ed\'eric MOTTE
\footnote{Amazon Science, \sf medericmotte at gmail.com 
}
\qquad\quad
Huy\^en PHAM
\footnote{LPSM, Universit\'e Paris Cit\'e, and CREST-ENSAE, \sf pham at lpsm.paris
This work was partially supported by the BNP-PAR Chair ``Futures of Quantitative Finance", and the 
Chair Finance \& Sustainable Development / the FiME Lab (Institut Europlace de Finance)
}
}

\begin{document}

\maketitle

\begin{abstract}
We study binary opinion formation in a large population where individuals are influenced by the opinions of other individuals. The population is characterised by the existence of (i) communities where individuals share some similar features, (ii) opinion leaders that may trigger unpredictable opinion shifts in the short term (iii) some degree of incomplete information in the observation of the individual or public opinion processes. In this setting, we study three different approximate mechanisms: common sampling approximation, independent sampling approximation, and, what will be our main focus in this paper, McKean-Vlasov (or mean-field) approximation. We show that all three approximations perform well in terms of different metrics that we introduce for measuring population level and individual level errors. In the presence of a common noise represented by the major influencers opinions processes, and despite the absence of  idiosyncratic noises,
we derive a propagation of chaos type result. For the particular case of a linear model and particular specifications of the major influencers opinion dynamics, 
we provide additional analysis, including long term behavior and fluctuations of the public opinion.
The theoretical results are complemented by some concrete examples and numerical analysis, illustrating the  formation of echo-chambers,   the propagation of chaos, and phenomena such as   snowball effect and social inertia.
\end{abstract}

\vspace{5mm}


\newpage

\section{Introduction}
Understanding how peer influence shapes the dynamics of opinions in a population is important in disciplines as diverse as economics, sociology, psychology, finance, medicine and innovation management.  
There is a substantial literature on this topic, to  be discussed below, that models peer influence as a contagion process.

However, the expansion of online social media is posing some new challenges to the existing models.  A clear  illustration are election outcomes that more and more often diverge from the predictions. Observers of the global political landscape have noted swift, drastic changes in opinion dynamics and polarisation, suggesting that internet in general, and social media in particular play a key role (\cite{BailAl18}, \cite{ZhuPetEni20}, \cite{Levy21}). But at the same time,  social media offers unprecedented possibilities for estimating peer influence through large scale data sets. Therefore, it appears that the theoretical tools need to be revisited, for a better integration of the features of the social networks that may be either observed or approximated reliably.  
What is special about online social influence, compared to traditional off-line influence? In this paper, we present a theoretical framework that integrates  aspects that we believe crucial, namely (i) the existence of communities that are randomly formed (hence do not match physical or socio-economic proximity), (ii) central role of opinion leaders that may create unpredictable opinion shifts  in the short term, (iii) some degree of incomplete information. 
 
 More precisely,  the considered aspects are as follows:

Firstly, individuals personalize their information access and tend to form communities that exhibit varying levels of interconnectivity, sometimes leading to reinforced beliefs within "echo chambers." These communities are often anonymous to external observers, and their membership is difficult to trace, distinguishing them from traditional groups that relied on geographical or socioeconomic proximity. To account for this aspect, we directly model influence at the community level, focusing on the distribution of individuals within communities rather than the underlying network based on individual identities.

Secondly, a new type of opinion leaders emerged in social media, with influential power over individuals (consumers, voters), transcending the communities.  Traditionally, opinion leaders used to be recognisable by their knowledge, expertise, or experience and their influence was more local (e.g., geographically). In online social media, influencers are  recognisable by the number of followers and centrality in the social network.
Therefore they are important determinants of rapid and sustained behaviour change.  In our model, a group of main influencers is introduced, having an exogenous  opinion dynamics that may be predictable or unpredictable, and driving the opinion of all the population. In practice, the shifts of the main influencers' opinions may be tactics to increase virality, or intent to initiate new fads. Opinions could also be interpreted as paid advertising, that can be diffused via diferent media. With this interpretation, these main influencers are the media, that can diffuse different advertising in time.

Thirdly, abundance of information means a high cost of processing it.  It is therefore not realistic to assume in the context of today's online influence that individuals are aware at all times of all reviews or opinions susceptible  to influence them. 
In this paper, we consider that individuals are potentially influenced by the whole population, but it is not possible for them to rely on a  full census of the population. As a substitute, individuals use an approximate mechanism to estimate the public opinion. We study three different approximate mechanisms: common sampling approximation, independent sampling approximation, and, what will be our main focus in this paper, McKean-Vlasov (or mean-field) approximation, the latter corresponding to our proposed concept of ``implicit influence''.  From the viewpoint of social science,  implicit influence may be interpreted in terms of social norms, or the feature of a population to have a set of (complex and stable) beliefs about how others will act in different social situations, and which guides individual behaviour in the absence of the observations of peer behaviour.

\paragraph{Related Literature.}
The study  of peer influence and opinion dynamics within networks has a long history, with typical applications focusing on social learning or adoption of behaviors and innovations. In a social learning context, it is assumed that there is a true state of the world -- initially unobserved --  and a network of individuals collectively
learns by aggregating information from their individual experiences. In adoption of innovation, there may not be an explicit unknown state of the world to be discovered, but instead agents utilities depend on coordination of consumption and hence more people adopt a product, more people may be incentivised to adopt in the future. Independently of the its application (learning or adoption),  there are two distinct modelling approaches. The first approach is to formulate a dynamic game and characterise its Bayesian equilibria (see the pioneering work by \cite{Ban92} and \cite{BikHirWel92}, or more recent work by \cite{SmiSor00}, \cite{BanFud04}, \cite{GalKar03}, Solan et al. \cite{SolRosVie09}, Acemoglu et al. \cite{AceDahLobOzd08}). Though theoretically attractive, the characterisation of such equilibria in complex networks is generally a non-tractable problem (see \cite{Hazal20}). Also, individuals in practice use somewhat coarser ways of aggregating information and forming beliefs than Bayesian updating.

These considerations have motivated the development of an alternative literature that explores belief formation based on  ``reasonable rules of thumb''. Our paper aligns with this approach.  The basic idea draws inspiration from  DeGroot's framework \cite{DeGroot74}, where individuals repeatedly average the opinions of others in the population. This can be synthesised as follows: when individuals gain access to new information, they aggregate it using an exogenously specified rule (usually simple to state and not time-varying)  to update their
current opinion. See for instance Ellison and Fudenberg \cite{EllFud93}, \cite{EllFud95}, Bala and Goyal \cite{BalGol98}, DeMarzo, Vayanos and Zwiebel \cite{DeMVayZwi03}, Golub and Jackson \cite{GolJac10},  Mossel et al. \cite{mosslytam14}, \cite{MossSlyTam15}, Popescu and Vaidya \cite{PopeVaid22}. In the particular case where  the opinion process is binary (i.e. adoption or not of an opinion by an individual),  the model is often of  ``threshold type'', i.e., an individual adopts once that sufficient neighbours have adopted (i.e., the individual's score exceeds a given threshold), see for instance \cite{Wat02}, \cite{Young09}, \cite{Young11}. 
The idea underlying the threshold models is as follows: even though the detailed mechanisms involved in binary decision problems can vary widely across specific problems, many decisions are inherently costly, requiring commitment of time or resources, hence the relevant decision function exhibits a threshold nature: agents display inertia in switching states, but once their personal threshold has been reached, the action of even a single neighbour can tip them from one state to another. Threshold models may also be obtained starting from  competitive games where coordination brings additional utility,  and adding additional noise in the player's payoffs \cite{AriBabRerPey2020}.
The interested reader may find a more extensive literature review of the models of social interaction in  the book by Jackson \cite{Jac08}, or the survey articles \cite{mostam17} and \cite{AceOzd}.

\paragraph{Our contributions.} 
We propose  an influence mechanism of threshold-type: opinion of each individual is binary and we do not explicitly assume that there is a true state of the world that needs to be discovered. In our model, a group (main influencers or fad) has an exogenously specified opinion dynamics, impacting potentially each individual in the larger population. 

 In our proposed setting, our main results are the following. 
 First, for  different approximate mechanisms of the population opinion dynamics (i.e., common sampling, independent sampling, and  McKean-Vlasov
approximation),  we provide simple estimations of how well these different approximates are close to the influence mechanism under complete information. This is very important in applications, whenever the costs of imple\-menting full information models are simply too high, in terms of data collection and/or processing.  Using an approximation mechanism for every individual, the resulting opinion dynamics follows a different path from the one obtained under perfect information. The question that arises naturally is whether the local (individual-level, or micro-scale) approxi\-mations compound to give large errors at the global level (population-level, or macro-scale). We will study the quality of an approximation via two natural metrics that we call the {\it local} and the {\it global} imperfectness of each mechanism.  Importantly, our finding is that all the  approximations that we have used, including the implicit influence mechanism (McKean-Vlasov) approximate well the complete information mechanism, in both metrics. 
From a mathematical viewpoint, such result contributes to a growing literature on mean-field inter\-actions with common noise, following the seminal work \cite{szni89}. Conditional propagation of chaos has been obtained in different models (e.g.  \cite{cogfla16}, \cite{cardellac16}, \cite{ernloclou21}). We refer to the monograph \cite{cardel} for a detailed treatment of the topic, together with an additional list of references, and also the recent survey in \cite{chadie22}.
In the current paper, although we  do not have idiosyncratic noises, but only a common noise represented by the major influencers, we can derive a conditional propagation of chaos type result.

Secondly, in the case of the McKean-Vlasov approximation and the linear macro-scale function case, we also study some long-term behavior of the public opinion processes. For this, we make two alternative assumptions for the dynamics of the main influencers opinion process: Markov chain assumption, versus deterministic, regular fluctuations. In the first case, we compute moments of the limiting distributions of the resulting public opinion, assuming a stationary distribution for the main influencers opinions. In the second case, we study the fluctuations size of the public opinion, in reaction to the changes in a main influencer's opinion process, that alternates with constant periodicity. 

Thirdly, we illustrate our results by examples and numerical simulations, revealing phenomena as such the eco-chambers effect, snowball effects or social inertia. The co-chamber effect is a polarisation effect: starting from a population where the opinion is i.i.d. distributed initially, the personal features of opinion formation can lead in the long-run to distinct communities adopting opposite opinions.  The snowball effect and social inertia refer to some  patterns of the  reactions in  the opinion processes to  changes of main influencer's opinions. More exactly, the snowball effect refers to the fact that social influence acts as an amplifier of the size of the fluctuations of the opinion processes of the population, while and social inertia is the property that social influence negatively impacts the speed of adjustment of a public opinion to changes in the main influencers' opinions.

Finally, we give an original interpretation of our mean-field model, in terms of a new concept of implicit social influence. 
Originated from statistical physics, mean-field models are usually viewed as an approximation of a physical system consisting of a large population of interacting  particles/agents. In our context, the physical system is a population, with a behavior (opinion formation) impacted by observations and beliefs. Hence, it is realistic to assume that individuals cannot observe the opinions of all the other individuals in the population and  they adapt behavior to the available observations. They use approximations, when the information is not complete. Hence, the approximations are not necessarily used by an external observer (as in statistical mechanics) but they impact the evolution of the physical system. The implicit social influence is an extreme case, where agents only observe the opinions of the main influencers and not of the remaining of the population. The corresponding approximation is the McKean-Vlasov.


 \paragraph{Outline of the paper.}
The model of the opinion dynamics in different information settings is presented in Section \ref{sec-Model}. Section \ref{sec-Estimate} introduces several types of error approximations, when information is incomplete, we then prove that for the different metrics the approximations work well. Section \ref{sec-Linear} presents in more detail the linear model, i.e., obtained when the macro-scale function is linear, and under some particular specifications of the dynamics of the main influencers. Theoretical results on long term behavior and fluctuations of the resulting opinion process are then described. We provide in Section \ref{sec:num} some  illustrative examples in a numerical study. Finally, all proofs are presented in Appendix \ref{sec-Proofs}.

\section{The model}\label{sec-Model}
A probability space $(\Omega,\mathcal F, \P)$ is fixed, on which all random variables will be defined; all relations involving random variables are considered in an a.s. sense. 

We consider a large population with $N$ interacting individuals, and denote by $X^N_n(t)$ the binary opinion, valued in $\{0,1\}$, of individual $n$ $\in$ $\llbracket 1,N\rrbracket$ at time $t$ $\in$ $\N$.  We may interpret opinion 0 as the status quo, while opinion 1 is the choice of a new technology, or innovation. With such an interpretation, the evolution of the stochastic process $\boX^{N} =(X^N_1(t),\ldots, X^N_N(t) )_{ t\in\N}$ reflects the adoption of the innovation by the population. As we are going to see below, the individuals may change their opinion multiple times, hence the setting is also suitable for modelling consumer choices, voting intentions for an election or spread of rumors in a population.

In addition to these $N$ individuals, we assume there are $M_0$ major influencers (or opinion leaders),  represented by a stochastic process  $\boX^0$ $=$ $(X^0_1(t),\ldots,X_{M_0}^0(t))_{t\in \N}$ valued in $\{0,1\}^{M_0}$, where $X_m^0(t)$ is the opinion of influencer $m$ $\in$ 
$\llbracket 1,M_0\rrbracket$ at time $t$ $\in$ $\N$. Major influencers are a very limited group of people, as compared to the size of the overall population, who are active in the interpersonal communication network, supplying infor\-mation, opinions, and suggestions, exerting personal influences on others. The term ``opinion leaders'' was coined by Lazarsfeld \cite{LazBerGau44}, when studying the effect of interpersonal commu\-nications  after the general election in 1940 in his US book ``The People’s Choice", while the term ``main influencers'' is more specific to online marketing. In our paper, we are considering the opinions of  the main influencers as exogenously given, and we introduce a mechanism by which  individuals $n\in \llbracket 1,N\rrbracket$ form their  own opinion.  This will allow us to study the speed of adjustment of the public opinion, in reaction to both predictable and non predictable variations of the process $\boX^0$.

The mechanism we propose for the opinion formation is based on two distinct parts: 
\begin{itemize}
\item[(1)] {\it Influence type}:  it will be captured by a set of personal features indicating how an individual is influenced by the others. These features are fixed for each individual, and will determine how the individual aggregates the information available to compute an own ``score''.

\item[(2)] {\it Information flow}, or information observable by an individual at each point in time; changes in the information flow will create the opinion dynamics (i.e., individuals update their ``scores''). We shall propose different models for the information flow. 
\end{itemize}
To fix  ideas, let us immediately present the full information setting where  the opinion processes  $\boX^N$ and $\boX^0$ are observable at time $t$ by all individuals, before turning to alternative information settings. With full information, we assume that at time $t$ any individual $n$ will aggregate  the observed opinions and compute  an own score, of the form:
\beqs
\mrS_n^N(t) &:=& \sum_{j\in \llbracket 1,N\rrbracket}\pi^N_{n,j}X_{j}(t)+ \sum_{m\in \llbracket 1,M^0\rrbracket} c_{n,m}^0 X_m^0(t)
\enqs
where the weights $(\pi^N_{i,j}), i,j\in\llbracket 1,N\rrbracket$ and $(c^0_{n,m}),  n\in\llbracket 1,N\rrbracket, m \in  \llbracket 1,M_0\rrbracket$, are random variables.  

In our framework, $(\pi^N_{i,j})$  reflect the existence of communities of opinion, where some fractions of the population share some similarities in the way they average opinions to form an own score.  More precisely, we shall assume that the population is partitioned in $K$ communities of influence, or classes represented by $\Kc$ $:=$ $\lb 1,K \rb$.  Let $(\kappa_i)_{i\in \llbracket 1,N\rrbracket}$ be random variables, such that $\kappa_i$ equals the class of individual $i$. We postulate that the influence that an individual $j$ exerts over another individual $i$ only depends on their respective communities. We  consider a matrix $\bolc$ $=$ $(c_{k,\ell})_{k,\ell \in \Kc}$ $\in$ $\R^{K\times K}$, 
representing the \textit{ inter-class influence matrix} and we assume 
\begin{equation}\label{eq-pi}
    \pi^N_{i,j}=\frac{1}{N}  c_{\kappa_i,\kappa_j}.
\end{equation}

In our setting, the influence matrix $\bolc$ is not random, but we assign each individual randomly to one of the $K$ classes, hence the resulting random weights $(\pi^N_{i,j})$.  Hence, individuals belonging to a class have similar updating rules, in the sense that their opinion will be influenced by the population in a similar way. 
The matrix $\bolc^0$ $=$ $(c^0_{n,m})_{n,m}$ may be random and the coefficient $c_{n,m}^0$ represents the impact of the main influencer $m$ on individual $n$. 

Consequently, the average opinion in communities will play an important role in the opinion formation and its time-dynamics. Indeed, we notice that with full information, the opinion score computed by individual $n$ at time $t$ has an alternative expression
\begin{align}\label{eq-score}
\mrS_n^N (t)  &:= 
\bolc_{\kappa_n}\cdot \bP^N(t)  + \bolc_n^0 \cdot \boX^0(t) 
\end{align}
where $\bolc_k$ $=$ $(c_{k,\ell})_{\ell\in\Kc}$ is the $k$-th row  of the matrix $\bolc$, 
$\bolc^0_n$ $=$ $(c^0_{n,m})_{m\in\llbracket 1,M^0\rrbracket}$ is the $n$-th row  of the matrix $\bolc^0$
and 
\beqs
\bP^N= (P^{N}_1(t), \dots, P^{N}_K(t))_{t\in \N}
\enqs
with 
\beqs
P_k^{N}(t)  = 
\frac{1}{N}\sum_{i=1}^N X^N_i(t){\bf 1}_{\kappa_i=k} 
\enqs
representing the proportion of opinion $1$ at time $t$ in each class $k$ $\in$ $\Kc$, expressed  over the whole population of $N$ individuals.  Notice that $\bP^N$ is valued in $\Sc_K$ $:=$ $\{ \bp=(p_k)_k \in [0,1]^K: \sum_{k=1}^K p_k \leq 1\}$. 

Given the above score computed by individual $n$,  we assume that she/he will then compare it with his/her  own  threshold $\mrs_n$, and will adopt opinion $1$ if and only if the score is above $\mrs_n$, i.e. $\mrS_n^N(t)$ $>$ $\mrs_n$. 
The interpretation of the threshold $\mrs_n$ is that individual $n$ initially has a natural inclination towards opinion $0$ or opinion $1$, which is represented by $\mrs_n$: indeed, the higher $\mrs_n$ is, the higher the score will have to be in order  to convince the individual to adopt opinion $1$.  

In addition to full information, several frameworks will be introduced below, where the population observes only partly peer opinions. But independently of the information available, there are some characteristics or features of the individuals that remain fixed. Hence the definition: 
\begin{Definition}\label{def:features} For any individual $n\in  \llbracket 1,N\rrbracket$, the features vector (or individual characteristics vector)
$$
\chi_n=(\xi_n,\kappa_n, \bolc^0_n=(c^0_{n,m})_{m\in\llbracket 1,M_0\rrbracket},\mrs_n)
$$ is a random vector  valued in $\mathcal X:= \{0,1\}\times \Kc\times\R^{M_0}\times\R$, where $\xi_n$ is the initial opinion of individual $n$, $\kappa_n$ denotes the class of individual $n$, 
$\bolc_n^0$ $=$ $(c^0_{n,m})_{m\in\llbracket 1,M_0\rrbracket}$ the influence coefficients (positive or negative) of the major influencers on individual $n$, and $\mrs_n$ represents the threshold of individual $n$.
\end{Definition}

\begin{Assumption}
The family of random vectors $(\chi_n)_n$ in Definition \ref{def:features} are supposed to be i.i.d. and also independent from the process $\boX^0$.
\end{Assumption}

Let us introduce some illustrative example of how different features of the population may be used to model phenomena of interest.

\begin{Example}[Echo-chambers]\label{ex:eco-chambers} Opinion polarisation is a process by which opposite opinions will tend to be adopted by distinct groups as time goes on. This phenomenon can be reproduced as follows.
Consider $(\xi_n)$ are i.i.d. Bernoulli with parameter $\epsilon$ (small), meaning that the opinion 1 is rare in the population originally at time 0. We consider a major influencer that spreads the opinion 1, by adopting an opinion process that is Bernoulli with parameter $(1-\epsilon)$. We consider two classes of agents and a corresponding inter-class matrix:
\[
c=\begin{pmatrix}
1/2 & 0\\
0 & 1
\end{pmatrix}.
\]
As there is only one main influencer, $\bolc^0_n\in\R$ and we assume that  these influence coefficients satisfy, conditionally on the class of the individual $n$:
\[
\begin{cases}
c^0_n=1/2 \text{ if }\kappa_n=1\\
-c^0_n\sim \text{Bernoulli}(\nu) \text{ if }\kappa_n=2
\end{cases}
\]
for some $\nu$ $\in$ $(0,1)$. 
The thresholds $(\mrs_n)$ are supposed uniform on [0,1]. In average, opinion 1 will become predominant in class 1 and it will tend to disappear in class 2, as will be discussed later on (see Example \ref{ex1}). 
\end{Example}

\begin{Example}[Fads with early adopters and late adopters]\label{ex:cycles}
We consider again a single main influencer, that adopts opinion 1 for a period $T$, then opinion $0$ for the following $T$ units of time, then reverts to opinion 1 and so on. This behavior may be interpreted as a fashion adoption that changes  cyclically. There are two classes of agents (early adopters, with $\kappa_n=1$ and late adopters, with $\kappa_n=2$) and initially, all the population has opinion 0. We assume the following influence coefficients:
\[
c=\begin{pmatrix}
1/2 & 0 \\
1/3 & 1/3
\end{pmatrix}
\]
    \[
\begin{cases}
c^0_n=1/2 \text{ if }\kappa_n=1\\
c^0_n=1/3 \text{ if }\kappa_n=2.
\end{cases}
\]
\end{Example}
The thresholds $(\mrs_n)$ are supposed uniform on [0,1]. 
 
 \begin{Remark}
     As mentioned in the introduction, it is common in some literature to interpret the scores as an average of neighbor opinions, hence take the sum of all weights equal to 1. In our setting, as we have random weights, we can implement this idea by setting 
     $\sum_{j=1}^N\pi_{i,j}^N+\sum_{m=1}^{M^0}c^0_{i,m} $
     to be close to 1 as $N$ large. This can be achieved as follows. Let $\mu=(\mu_k)$ be the probability mass function of $\kappa_n$. Then, we impose
     \[
     \sum_{k=1}^K\mu_k c_{\kappa_i,k}+\sum_{m=1}^{M^0}c^0_{i,m}=1
     \]which would lead to the desired property by \eqref{eq-pi} and the strong law of large numbers. To fulfill the constraints, the random weights $\{c^0_{i,m},m\in \llbracket 1,M^0\rrbracket \}$ may be chosen $\sigma(\kappa_i)$ measurable, which would mean that the main influencers impact is at the class level. However, our framework is more general, and does not require this specification.
 \end{Remark}

\vspace{3mm}

We now propose several frameworks for the information available and the corresponding scores and opinions.

\subsection{Full information }

We start with the description of   the population's opinions dynamic in the  case where each individual has complete information, i.e. can access to the state of  the whole 
population at each time.  In this case,  the information available to any individual at time $t$ is given by:
$\sigma (\boX^N(s), \boX^0(s), s<t)\bigvee_{n\in \llbracket 1,N\rrbracket} \sigma(\chi_n)$. 
 We assume in this case that the corresponding score of individual $n\in \llbracket 1,N\rrbracket$ at time $t$ is as  in \eqref{eq-score} with the following opinion dynamics: 
\begin{align}\label{eq-simple-dyn}
\begin{cases}
X^{N}_n(0) &= \; \xi_n,  \\
X^{N}_n(t+1) &=  \;   {\bf 1}_{\mrS_n^N(t) > \mrs_n}, \quad  \; t \in \N, 
\end{cases}
\end{align}
where we recall that $(\xi_n)_{n\in \llbracket 1,N\rrbracket}$ are i.i.d. random Bernoulli variables representing the initial opinions in the population, with the interpretation that people's initial opinions are independently generated by a same random ``genetic'' mechanism. 

\subsection{Complete but asynchronous observations}
Although the above dynamics is simple and natural, it may lack of realism because it assumes that all the population is influenced simultaneously on a discrete time grid. To make it more realistic, a basic  extension would be to  
consider that at each day $t$, an individual  is influenced 
only with probability $h$ $\in$ $(0,1]$. On modern media like social networks and internet in general, individuals generally connect to the web with a much smaller degree of synchronization: this can be modeled by considering a very small $h\simeq 0$, for which the event that two individuals connect at the same time will have a very small probability. 

In order to deal with the above issues, we shall generalize the dynamic \eqref{eq-simple-dyn} as follows. We introduce the notion of  {\it influence event}: we consider a family $\{B_n(t), t\in\N\},n\in \lb 1,N \rb$ of i.i.d. Bernoulli processes with parameter $h$, where the event  ``$B_n(t)=1$'' (resp. ``$B_n(t)=0$'') means that at time $t$, individual $n$ connects (resp. does not connect) to the web (or watches (resp. does not watch) television).  In other words,  each individual $n$ is influenced by the population at 
successive times $0=\tau_n^0 < \tau_n^1 <  \ldots < \tau_{n}^j < \ldots$, and  $\{\tau_{n}^{j+1}-\tau_n^{j}:  j  \in \N,  n \in \llbracket 1, N\rrbracket \}$ is a family of  i.i.d. inter-arrival random times  with marginal distribution $Geom(h)$, the geometric distribution of parameter $h$ $\in$ $(0,1]$. 
We assume that the family $\{B_n(t): n\in \lb 1,N \rb, t\in\N\}$ is independent of the family $(\chi_n)_{n\in \llbracket, 1,N\rrbracket}$ and of the process $\boX^0$.

We assume that with complete, asynchronous observations, new information is available to agent $n$ only at its influence times, and then he/she will be able to observe the opinions of the whole population.
Consequently,  the score computed by individual $n$ equals $\mrS^N_n(t)$,  defined in \eqref{eq-score} whenever $B_n(t)=1$ (i.e., at any influence time), and the individual is not updating the opinion when $B_n(t)=0$ (hence is not needing the score). We now define  the generalized version of the opinion dynamic with complete information, as the stochastic process  $\bX^{N}= (X^N_n(t),n \in \llbracket 1,N\rrbracket )_{t \in \N}$,  given by
\begin{align} 
\begin{cases}
X^{N}_n(0)   = \;\xi_n, \\
X^{N}_n(t+1)  =\; (1-B_n(t))X^{N}_n(t)+B_n(t){\bf 1}_{\mrS_n^N(t) > \mrs_n},  \quad t \in \N.  \label{eq-general-dyn}
\end{cases}
\end{align}
where $\mrS^N_n(t)$ is defined in \eqref{eq-score}. Notice that one retrieves the synchronized dynamic \eqref{eq-simple-dyn} when $h=1$.


\subsection{Incomplete information and approximations} 
The assumption that individuals had a complete access to the opinions of all the population may also be criticised for various reasons;  (i) 
{\it the size of the population}: in large populations, it is too costly to access and process the opinions of  the whole population; (ii) {\it the frequency}: in situations where individuals may potentially change opinions multiple times, it is simply too costly to access the opinions of the population at frequent repeated times; (iii) {\it privacy issues}: some individuals might not be willing to openly share their opinions. 

Therefore, we propose a more realistic setup, that does not require individuals to survey the whole population to form a public opinion. More exactly, we propose three alternative methods for approximating the vector process $\bP^N$, i.e., the average opinion in communities, our assumption being that each individual $n$ will use an own approximation $\tilde \bP_n$, that will be a stochastic process adapted to an own information flow of agent $n$. By contrast, the opinion of the main influencers is always assumed to be observed by all individuals.

An approximate model of the public opinion process is formalized by an opinion process of the population $\boX$ $=$ 
$(\tilde X_1(t),\ldots,\tilde X_N(t))_{t\in\N}$ governed by 
\begin{equation} \label{dynXtildep}
\begin{cases}
\tilde X^{}_n(0)   = \;\xi_n,  \quad  n =1,\ldots,N, \\
\tilde X^{}_n(t+1)  =\; (1-B_n(t)) \tilde X^{}_n(t) + B_n(t)   {\bf 1}_{\tilde\mrS_n(t) > \mrs_n}, \quad  \;  t \in \N, n =1,\ldots,N,   
\end{cases}
\end{equation}
with approximate score 
\beqs
\tilde\mrS_n(t) &:=& \bolc_{\kappa_n}\cdot \tilde \bP_n(t) + \bolc_n^0\cdot \boX^0(t), 
\enqs
for some process $\tilde\bP_n$ $=$ $(\tilde{P}_{n,1}(t),\ldots, \tilde{P}_{n,K}(t))_{t\in\N}$ valued in $\Sc_K$, and adapted to the information  flow of individual $n$.

We shall focus on  three types of public opinion approximation mechanisms, depending upon what information individuals are assumed to  access, namely:
\begin{itemize}
\item[(1)] {\it Common $M$-sample groups approximation}. In this case, at any time $t$ all individuals that update their score are using  the same approximations, obtained from surveying $M$ individuals chosen uniformly at random from the entire population. This corresponds to 
\begin{align} \label{independent} 
\tilde P_{n,k}(t)  &= \;  \frac{1}{M} \sum_{i \in I(t)} \tilde X_i^{}(t){\bf 1}_{\kappa_i=k}, \quad k \in \lb  1, K \rb\;,\;n \in \lb  1, N \rb, 
\end{align} 
where  $(I(t))_{t\in\N}$ is a sequence of random subsets of $M$ $\leq$ $N$  indices uniformly sampled in the set $\{1,...,N\}$, without replacement. We notice that indeed $\tilde P_{n,k}(t) $ are not depending on $n$ (at a given time, two individuals use the same approximations).

\item[(2)] {\it Independent $M$-sample groups approximation}.  This corresponds to the case where each individual makes an own survey of size $M$, independently of other individuals  and use it as an approximation.
\begin{align} \label{sample} 
\tilde P_{n,k}(t) &= \;   \frac{1}{M} \sum_{i \in I_n(t)} \tilde X_i^{}(t){\bf 1}_{\kappa_i=k}, \quad k \in \lb  1, K \rb\;,\ n \in \lb  1, N \rb, 
\end{align} 
where $\{I_n(t): n \in \llbracket 1,N\rrbracket, t \in \N\}$ is a family of i.i.d.  random subsets of $M$ indices uniformly sampled in $\{1,...,N\}$, without replacement. 

\item[(3)] {\it A mean-field (or McKean-Vlasov) influence approximation}. Here, all individuals observe the opinion of the main influencers only, hence use identical approximations, as follows:  for $ k \in \lb  1, K \rb\;,\;  \tilde P_{n,k}= P^{\mbox{{\tiny{MKV}}}}_k$, for all $n \in \lb  1, N \rb$, where:
\begin{align} \label{tildeMKV}
 P^{\mbox{{\tiny{MKV}}}}_{k}(t) &=   \E \big[  \tilde X^{}_{n}(t) {\bf 1}_{\kappa_n=k}\mid \boX^0(s),s\leq t \big], \quad k \in \lb  1, K \rb\;,\;  n \in \lb  1, N \rb. 
\end{align}
By plugging \eqref{tildeMKV} into \eqref{dynXtildep}, it is indeed clear that $\tilde X_n^{}$, $n$ $\in$ $\llbracket 1,N\rrbracket$ are identically distributed, hence $\tilde P_{n,k}$ does not depend on $n$, see more details below. 
We denote $\bP^{\mbox{\tiny{MKV}}}$ $=$ $(P_1^{\mbox{\tiny{MKV}}}(t),\ldots,P_K^{\mbox{\tiny{MKV}}}(t))_{t\in\N}$. 
\end{itemize}

Let us discuss in more details the public opinion approximation mechanisms. Consider a fixed time $t\in\N$ where  some individuals (i.e., with indices belonging to $\{n\in \lb{1,N}\rb: B_n(t)=1\}$) want to update their score. Without the possibility to observe the opinion of the whole population, these individuals are assumed to rely on one of the three approximations introduced above.

In the common sample groups opinion dynamic, at each time $t$, a group of $M$ individuals, with indices in $I(t)$, is randomly selected from the general population  and their opinions and class types are observed. All individuals updating their score are assumed to use 
the quantity $(\frac{1}{M} \sum_{i \in I(t)} \tilde X_i^{}(t){\bf 1}_{\kappa_i=k})_{k\in\Kc}$  instead of $(\frac{1}{N}\sum_{i=1}^N X^N_i(t){\bf 1}_{\kappa_i=k})_{k\in\Kc}$ within in the score calculation. This is thus a natural representation of situation where regular surveys are published or updated. It emulates standard surveys communicated on TV show, but also the likes and dislikes, for the online reviews or social media. 

In the independent sample groups opinion dynamic, at any time, each individual samples a group of $M$ individuals at any time of influence. The difference with the previous case is that not only the sampled groups differ   from one influence time to the next, but also from one individual to another. 
This ``resampling'' of the sample group can have various causes.  One of them is when there are no officially organised surveys, or when the organised surveys are not precise enough, e.g., as mentioned above, when there are multiple classes and the official surveys only provide the global proportion of opinion $1$. In this case, people can make their own ``surveys'' by consulting forums and comments, and they can themselves class the authors of comments into groups (again, young, old, etc).

Finally, the McKean-Vlasov dynamic corresponds to the concept of implicit influence, that is our main focus in this paper. Implicit influence applies in a situation where people cannot access peer opinions, even of a sample group, but only the opinion of the main influencers at the influence times.  Then, people will unconsciously  theoretically model the population to predict its evolution as time goes. We assume that the way the individual models the population is via a mean-field approximation (McKean-Vlasov). An alternative use of the McKean-Vlasov approximation is in the following situation. There is an external observer that has interest in knowing the public opinion, which is however unaccessible. The external observer can only access the main influencer's opinions and has knowledge of the distribution functions of the personal  features in the population.  The observer also knows that individuals only update their information occasionally, that is with probability $h$ at a point in time, independently one from another. In any of the two settings, the approximation errors are useful to understand.

In the McKean-Vlasov dynamic, the only information that individuals access in order to form an opinion are the opinions of the major influencers, i.e. $\boX^0$. Hence, the resulting opinion processes $(\tilde X^{}_n(t)), n\in  \llbracket 1,N\rrbracket$ will be independent conditionally to $\boX^0$.  Indeed, recall that the only other ``information'' that each individual $n$ has is her/his own features vector $\chi_n$,  and $(\chi_n)$ are assumed to be i.i.d. and to also remain i.i.d conditionally to $\boX^0$. 
Therefore, in \eqref{tildeMKV}, one has $\tilde \bP_n(t)$ $=$ $\tilde \bP_1(t)$, for all $n$ $=$ $1,\ldots,N$.

Intuitively, by the (conditional) law of large numbers, for all $n \in \llbracket 1,N\rrbracket, t\in \N$, it is expected that the public opinion $\frac{1}{N}\sum_{i=1}^N X^N_i(t) {\bf 1}_{\kappa_i=k}$ will be close to $\E[\tilde X^{}_{1}(t) {\bf 1}_{k_1=k} \mid \boX^0(s),s\leq t]$, for $k$ $\in$ $\Kc$, 
which explains intuitively why  individuals will naturally use this quantity to approximate the public opinion. The rigorous justification and quantification of the approximate error are given in the next section. 
We stress that the McKean-Vlasov approximation provides a dynamic model where social influence plays an important role while individuals do not even access the opinions of other people. This is why we refer to this phenomenon as \textit{ implicit social influence}.

\section{Estimate of approximation errors}\label{sec-Estimate}

We shall measure the quality of such approximation from a {\it local} and {\it global} perspective. 
For an approximate public opinion model $(\tilde\boX,\tilde \bP)$, we define the local error function (or micro-scale error function) at time $T$  by 
\beqs
\Ic^N_{\ell oc}(\tilde \bP,T) &=&  \sup_{n\in \llbracket 1,N\rrbracket}\E\Big[ \sum_{k\in \Kc} \big\vert \tilde{P}_{n,k}(T) - \frac{1}{N}\sum_{i=1}^N \tilde X^{}_i(T) {\bf 1}_{\kappa_i=k} \big\vert \Big],
\enqs
and its global (or macro-scale) error function at time $T$ $\in$ $\R_+$ by 
\beqs 
\Ic^N(\tilde{\bP},T) &=& \E\Big[\frac{1}{N}\sum_{i= 1}^N \big\vert \tilde X^{}_i(T) - X^{N}_i(T)\big\vert \Big]. 
\enqs 
The interpretation of these error functions is the following. $\Ic_{\ell oc}^N(\tilde{\bP},T)$ measures the expected distances between the approximation $\tilde{\bP}_n(t)$ of public opinion used by any individual $n$, 
and the real public opinion $\Big(\frac{1}{N}\sum_{i=1}^N \tilde X^{}_i(t) {\bf 1}_{\kappa_i=k}\Big)_{k\in\Kc}$, provided that the approximation used by all individuals is $\tilde \bP$, if individual is influenced at time $t$, i.e. $B_n(t)=1$. It thus measures how good the estimations made by individuals are in expectation, if their goal is to approximate the public opinion. The purpose of such measure is to give an upper bound to the maximal expected estimation error that individuals should be willing to make in order to accept using this approximation.

On the other hand, $\Ic^N(\tilde \bP,T)$ measures at time $T$ the proportion of individual who have a different opinion than what they would have if the population had followed the dynamic with complete information from the beginning. The purpose of $\Ic^N$ is to help understanding the longer term impact (i.e., that lasts for more than one period)  on the population's dynamic of the fact that individuals use an approximate influence mechanism, instead of the perfect information one. For instance, how much can the result of presidential elections that is to be held in a few months be affected by the use of a given approximation? $\Ic^N$ thus allows us to  study the robustness of the population's dynamic to various public opinion approximation mechanisms. Our main results state an estimate of the error due to the approximate mechanism for opinion dynamic.

\subsection{Local error}

The estimate for the local approximation error is independent of the time $T$ considered.

\begin{Theorem}[Local error] \label{theo-loc}
 Let $\tilde \bP$ be any of the following: (i) the common $M$-sample group approximation, (ii) the independent $M$-sample group approximation. Then:
    \beqs 
    \Ic_{\ell oc}^N(\tilde \bP,T) &\leq&   \frac{\sqrt{K}}{2}   \frac{1}{\sqrt{M}}. 
    \enqs 
     Let  $\tilde \bP$ be the McKean-Vlasov approximation. Then:
    \beqs 
    \Ic_{\ell oc}^N(\tilde \bP,T) &\leq&  \frac{\sqrt{K}}{2}   \frac{1}{\sqrt{N}}. 
    \enqs 
\end{Theorem}
For the proof of this result, see Subsection \ref{proof:theo-loc}.

\subsection{Mean influence  function and global error}

An important object in our next results will be played by  the so-called mean influence 
function $\boldsymbol \phi=(\phi_1,\dots \phi_K)$. 
More exactly, in our definition, the function  $\boldsymbol \phi$ $:$ $\Sc_K\times \{0,1\}^{M_0} \rightarrow [0,1]^K$, has components $\phi_k$ defined by
\beqs 
\phi_k(\bp, \bolx^0) &=& \P\big[ \bolc_{k}\cdot \bp+\bolc^0_n\cdot \bolx^0>\mrs_n \; , \;  \kappa_n=k\big],  \quad \forall \bp\in \Sc_K, \; \bolx^0\in \{0,1\}^{M_0},\; k\in\Kc. 
\enqs 
 
This function approximates for the full information case, the public opinion of the communities that would prevail at time $t+1$,  assuming that at the time $t$  the major influencers have opinions $\bolx^0=(x^0_1,\dots x^0_{M^0})$ and  the mean public opinion in communities is equal to $\bp=(p_1,\dots,p_K)$.  
 Indeed, conditionally to $\{\boX^0(t)=\bolx^0,\bP^N(t)=\bp\}$,
we have that the full information score in \eqref{eq-score} writes:
\begin{align*}
 \frac{1}{N}\sum_{i=1}^N c_{\kappa_n,\kappa_i}X^N_i(t)+ \bolc^0_n\cdot \boX^0(t) 
 & = \; \sum_{\ell\in \Kc} c_{\kappa_n,\ell}\cdot\Big(\frac{1}{N}\sum_{i=1}^N X^N_i(t){\bf 1}_{\kappa_i=\ell}\Big) + \bolc^0_n\cdot \boX^0(t) \\
 & = \;   \sum_{\ell\in \Kc} c_{\kappa_n,\ell} \;p_\ell + \bolc^0_n\cdot \bolx^0 \\
 & = \;   \bolc_{\kappa_n} \cdot \bp + \bolc^0_n\cdot \bolx^0.  
\end{align*}
It follows that any individual $n$ in the large population is adopting opinion $1$ at time $t$ if and only if
\begin{align*}
 &  \bolc_{\kappa_n} \cdot \bp + \bolc^0_n\cdot \bolx^0  \; > \; \mrs_n, 
\end{align*}
and thus the empirical average opinion in the community $k$ would be equal to the random variable
\begin{align}\label{eq-lfgn}
\Phi^N_k(\bp,\bolx^0) & := \;  \frac{1}{N}\sum_{n=1}^N{\bf 1}_{\bolc_{k}\cdot  \bp + \bolc^0_n\cdot \bolx^0>\mrs_n}{\bf 1}_{\kappa_n=k}.
\end{align}
We shall call the function $\boPhi^N$ $=$ $(\Phi^N_1,\dots, \Phi^N_K)$ the empirical average opinion (or the empirical macro-scale function), and observe that: 
\begin{align}
\Phi^N_k(\bp,\bolx^0) & \simeq \;  \P\big[ \bolc_{k}\cdot \bp + \bolc^0_n\cdot \bolx^0>\mrs_n,{\bf 1}_{\kappa_n=k}\big]  
\; = \; \phi_k(\bp,\bolx^0), 
\end{align}
where ``$\simeq$'' comes from the law of large numbers for $N$ large.  
 
 We make a  Lipschitz assumption on the   macro-scale influence function: there exists some positive constant $K_\phi$ such that 
\beqs 
\hspace{-1cm} {\bf (H_\phi)} \quad\quad\quad  \vert\boldsymbol  \phi(\bp,\bolx^0)-\boldsymbol \phi(\bq,\bolx^0)\vert_{_1}  \; \leq \;  K_\phi \vert \bp-\bq\vert_{_1}, \quad \forall \bp,\bq\in \Sc_K, \;  \bolx^0 \in\{0,1\}^{M_0}.
\enqs
Here $|.|_{_1}$ is the $\ell_1$-norm on $\R^K$. Examples of explicit  macro-scale influence functions will be given in the next sections.

Denoting by $\|\boPhi^N-\boldsymbol \phi\|$ $=$ $\sup_{\bp\in \Sc_K,\bolx^0 \in \{0,1\}^{M_0}} | \boPhi^N(\bp,\bolx^0) -\boldsymbol \phi(\bp,\bolx^0)|_{_1}$, we recall by Glivenko-Cantelli theorem that $\E[\Vert \Phi^N-\phi\Vert]$ converges to zero as $N$ goes to infinity.  
Moreover, we know from Dvoretzky–Kiefer–Wolfowitz inequality (see e.g. \cite{Dvoretzky:1956aa} and \cite{Massart:1990aa}) that  for all $\eps$ $>$ $0$:  
\beqs 
\P[ \Vert \boPhi^N-\boldsymbol \phi\Vert>\eps] &<& 2 K e^{-2N\eps^2},
\enqs 
which implies that for all $\alpha>0$, there exists a constant $C_\alpha$ such that
\beqs 
\E[\Vert \boPhi^N-\boldsymbol \phi\Vert] &<&  \frac{C_\alpha}{N^{\frac{1}{2}-\alpha}}.
\enqs 

\begin{Theorem}[Global error] \label{theo-glob}
\label{theo-disc}
Assume that ${\bf (H_\phi)}$ holds true.  
\begin{enumerate}[(i)]
    \item Let $\tilde \bP$ be the common or the independent $M$-sample group approximation, then
    \beqs 
    \Ic^N(\tilde \bP,T) &\leq& 
    \frac{K_\phi^T - 1}{K_\phi -1} 
    \Oc\Big(\frac{1}{\sqrt{M}}+ \E[\Vert \boPhi^N-\boldsymbol \phi\Vert ]\Big).
    \enqs 
    \item Let $\tilde \bP$ be the McKean-Vlasov approximation, then 
    \beqs 
    \Ic^N(\tilde \bP,T) &\leq&  
     \frac{K_\phi^T - 1}{K_\phi -1} 
    \Oc\Big(\frac{1}{\sqrt{N}}+ \E[\Vert \boPhi^N-\boldsymbol \phi\Vert ]\Big). 
    \enqs
\end{enumerate}
Here, we use the convention that $\frac{\ell^T-1}{\ell-1}$ $=$ $T$ when $\ell$ $=$ $1$, and the Big Landau terms $O(.)$ 
depend only on $K_\phi$, and $K$.
\end{Theorem}

The results in the above theorem are proven in Subsection \ref{proof:theo-glob}; Theorem \ref{theo-loc} and \ref{theo-glob} show that by using the approximate mechanisms  the population globally and the individuals, locally, are not compounding large errors through time. 
 In other words, the dynamics of the opinions in the population will not be very far apart from the perfect information dynamics. 
 A natural conclusion from this is that, provided that in the absence of surveys, individuals would use a McKean-Vlasov approximation instead, then the publication or not of regular surveys before an election should not impact a lot its result, as long as $N$ and $M$ are large. Another remarkable conclusion from Theorem \ref{theo-disc} is that if $\boldsymbol{\phi}$ is strictly  contracting in $\bp$, i.e. $K_\phi < 1$, the global error still converges to zero  for large $T$. 
Otherwise, if $K_\phi$ $\geq$ $1$, 
the global error can diverge  exponentially (for $K_\phi$ $>$ $1$) or linearly (for $K_\phi$ $=$ $1$), when $T\rightarrow +\infty$, 
see Figure \ref{fig:div}.

\begin{Remark}
Suppose that instead of the McKean-Vlasov approximation, the population uses another public opinion approximation.  
\begin{enumerate}
\item If the approximation used is the real public opinion (under complete information hypothesis) $\tilde\bP$ $=$ $\bP^N$,  we see that for all $T$ $>$ $0$,  
\beqs 
\E[\vert \bP^N(T)- \bP^{\mbox{{\tiny{MKV}}}}(T)\vert_{_1} ]&\leq& \sum_{k=1}^K \E\Big[\frac{1}{N}\sum_{i=1}^N \vert  X^N_i(T)
-\tilde X^{}_i(T)\vert 1_{\kappa_i =k} \Big] \\
& & \quad + \; \sum_{k=1}^K  \E\Big[ \big\vert \frac{1}{N}\sum_{i=1}^N  \tilde X^{}_i(T) 1_{\kappa_i =k} - P_k^{\mbox{{\tiny{MKV}}}}(T)\big\vert  \Big] \\
&\leq&  \Ic^N(\bP^{\mbox{{\tiny{MKV}}}}, T)+\Ic^N_{\ell oc}(\bP^{\mbox{{\tiny{MKV}}}},T). 
\enqs
Hence, when $K_\phi$ $<$ $1$ (see examples in the next section), this implies with Theorems \ref{theo-loc} and 
\ref{theo-disc} that $\E[\vert \bP^N(T)- \bP^{\mbox{{\tiny{MKV}}}}(T)\vert_{_1} ]$ is of order $\Oc\Big(\frac{1}{\sqrt{N}}+ \E[\Vert \boPhi^N-\boldsymbol\phi\Vert ]\Big)$ uniformly in $T$. 
\item If  the approximation used is the common or independent $M$-sample group surveys $\tilde\bP$ in \eqref{independent} or \eqref{sample},   we see   that for all $T$ $>$ $0$, $n$ $\in$ $\N$, 
\beqs 
\E\Big[  \big\vert\tilde\bP_n(T)- \bP^{\mbox{{\tiny{MKV}}}}(T)\big\vert_{_1}  \Big]
&\leq &  
\Ic^N_{\ell oc}(\tilde\bP,T) + 
 \Ic^N(\tilde \bP,T)+  \Ic^N(\bP^{\mbox{{\tiny{MKV}}}},T)+ \Ic^N_{\ell oc} (\bP^{\mbox{{\tiny{MKV}}}},T), 
\enqs 
and again, when $K_\phi$ $<$ $1$, we see from Theorems \ref{theo-loc} and 
\ref{theo-disc} that $\E\Big[  \big\vert\tilde\bP_n(T)- \bP^{\mbox{{\tiny{MKV}}}}(T)\big\vert_{_1}  \Big]$ is of order 
$\Oc\Big(  \frac{1}{\sqrt{M}} + \frac{1}{\sqrt{N}}+ \E[\Vert \boPhi^N-\boldsymbol\phi\Vert ]\Big)$ uniformly in $T$.
\end{enumerate}
This means that for $N,M$ large, all the considered dynamics will closely follow the McKean-Vlasov macro-scale dynamic $(P^{\mbox{{\tiny{MKV}}}}(t))_{t\in\N}$ described above, with an explicit rate of convergence.   
\end{Remark}

\subsection{McKean-Vlasov approximate opinion dynamics}

We now provide an analytical calculation of  the McKean-Vlasov  approximate  public opinion $\tilde\bP$ $=$ $\bP^{\mbox{\tiny{MKV}}}=(P_1^{\mbox{\tiny{MKV}}}(t),\ldots,P_K^{\mbox{\tiny{MKV}}}(t))_{t\in\N}$ at any time $t$. We recall that  all agents use the same  approximation processes in this case, that is $\tilde P_{n,k}= P^{\mbox{\tiny{MKV}}}_k$ for all $n$ and for $k\in \mathcal K$.

\begin{Proposition} \label{propP}
 If $\tilde \bP$ $=$  $\bP^{\mbox{{\tiny{MKV}}}}$ is the McKean-Vlasov public opinion approximation in \eqref{dynXtildep}, then it satisfies the relation
\beqs 
\bP^{\mbox{{\tiny{MKV}}}}(t+1) &=& (1-h)\bP^{\mbox{{\tiny{MKV}}}}(t) + h \boldsymbol \phi(\bP^{\mbox{{\tiny{MKV}}}}(t), \boX^0(t)),\quad \forall t\in\N. 
\enqs 
\end{Proposition}
{\bf Proof.} We have that $P_k^{\mbox{{\tiny{MKV}}}}(0)$ $=$  $\P[ \xi_n =1, \kappa_n=k\mid \boX^0(0) ]$ $=$ $\P[\xi_n =1, \kappa_n=k]$, that is independent from $n$,  as the features vectors are i.i.d. by assumption. Denote by $\Fc^0(t)$ $=$ $\sigma(\boX^0(s), s \leq t)$, $t$ $\in$ $\N$, the filtration generated by $\boX^0$. 
 From the definition of $\tilde \bP$ in \eqref{tildeMKV}, and $X^{\tilde \bP}$ in \eqref{dynXtildep},  we have
 \beqs 
P_k^{\mbox{{\tiny{MKV}}}}(t+1)&=& \P[ \tilde X^{}_n(t+1) =1, \kappa_n=k\mid \Fc^0(t+1) ] \\\; 
&=&\P\Big[ (1-B_n(t))\tilde X^{}_n(t)+B_n(t){\bf 1}_{\tilde \mrS_n(t) > \mrs_n} = 1, \kappa_n=k\mid \Fc^0(t+1) \Big] \\
&=& (1-h)\P\big[ \tilde X^{}_n(t) = 1, \kappa_n=k\mid \Fc_t^0 \big]  +  h\P\big[ \tilde\mrS_n(t)>s_n, \kappa_n=k\mid \boX^0(t)\big] \\
&=&(1-h)P_k^{\mbox{{\tiny{MKV}}}}(t)+h\phi_k(\bP^{\mbox{{\tiny{MKV}}}}(t), \boX^0(t)), 
\enqs 
which concludes the proof.
\ep

\vspace{3mm}

The above proposition provides a simple procedure for computing   the McKean-Vlasov approximation of the opinion dynamic. 
Thus, whenever individuals cannot retrieve any  information about the average opinion in the population at the influence times, each individual may calculate an inductively defined sequence as in Proposition \ref{propP}.  
This simplicity  for the computation of $P^{\mbox{{\tiny{MKV}}}}$ is an argument in favour of the McKean-Vlasov approximation mechanism. In the sequel, we provide some examples  of explicit calculations for $P^{\mbox{{\tiny{MKV}}}}$, and applications.


\section{Applications with linear mean-field approximation}\label{sec-Linear}

In this section, we study in more details a linear dynamics for the mean-field model, and emphasise its interest in applications. In particular, we are interested in two aspects (i) the limiting distribution of the population opinion and (ii) the fluctuations of the opinion process, depending on those of the main influencers.

For simplicity, we consider the case of perfect synchronized times, i.e. $h=1$, and assume that the impact of the main influencers on each individual is only via its class, i.e. $c_{n,m}^0$ $=$ $c^0_{\kappa_n,m}$ (by misuse of notation). 
We then denote by $\bolc^0$ $=$ $(c_{k,m}^0)_{k,m}$ the constant matrix in $\R^{K\times M_0}$, and 
$\bolc^0_k$ $=$ $(c_{k,m}^0)_m$ the vector in $\R^{M_0}$ assumed to be nonzero. 
The thresholds $(\mrs_n)_n$ are assumed uniformly distributed in $[0,1]$, and independent of 
the class features $(\kappa_n)_n$ distributed according to some distribution $\mu$ $=$ $(\mu_k)_k$ on $\Kc$. We denote by $M_K$ the diagonal $K\times K$-matrix with diagonal elements $\mu_k$. 
We also suppose that $c_{k,\ell},\bolc^0_{k,m}\geq 0$ for all $k,\ell\in \Kc$, $m\in \llbracket 1,M^0\rrbracket$ 
and $\sum_{\ell} c_{k\ell} + \sum_m c_{km}^0$ $\leq$ $1$, for all $k$ $\in$ $\Kc$, which implies that  
$\sum_{\ell} c_{k\ell}$ $<$ $1$ since $\bolc^0_k$ is nonzero. 

In this case, the macro-scale influence function is given by 
\beqs
\boldsymbol \phi(\bp,\bolx^0) &=&  \bC p + \boCc^0 \bolx^0, \quad p \in \Sc_K, \; \bolx^0 \in \{0,1\}^{M_0},
\enqs
where we set $\bC$ $:=$ $M_K \bolc$, and $\boCc^0$ $:=$ $M_K \bolc^0$,  
and clearly satisfies assumption ${\bf (H_\phi)}$ with $K_\phi$ $=$ $\max_{k,\ell}c_{k,\ell}$. Let us introduce 
\beqs
\bP^0(t) &:=&   \boCc^0  \boX^0(t),  \quad  t\in \mathbb N,
\enqs
that is   the weighted average opinion of the main influencers, re-scaled by the weight of each class.

The function $\bP^{\mbox{{\tiny{MKV}}}}$ defining the McKean-Vlasov approximate public opinion is then given from Proposition \ref{propP} by
\begin{align} \label{PMKVlinmulti}  
\bP^{\mbox{{\tiny{MKV}}}}(t+1) &= \;   \bC \bP^{\mbox{{\tiny{MKV}}}}(t) + \bP^0(t), 
\end{align} 
that leads to: 
\begin{align}\label{eq-attraction}
\bP^{\mbox{{\tiny{MKV}}}}(t) &=\;  \bC^t \bP^{\mbox{{\tiny{MKV}}}}(0)  + \sum_{s=0}^{t-1} 
\bC^s \bP^0(t-1-s).
\end{align}
Hence, the implicit McKean-Vlasov influence mechanism can be accessed by any individual in the population that is observing the trajectories of the major influencers. Therefore, the dynamics of the process $\bP^0$ is fundamental in approaching our questions of interest, such as convergence or fluctuations of the process $\bP^{\mbox{{\tiny{MKV}}}}$.
We may first look at what happens on the intervals where $\bP^0$ remains constant (assuming such a non empty interval exists). We introduce the sequence of stopping times $(\sigma^0_k)$ as: 
\begin{align}\label{eq-sigma0}
    \sigma^0_0&= \; 0, \quad \quad\sigma^0_{k+1} \; = \; \inf\{ t  >  \sigma^0_{k}: \bP^0(t) \not = \bP^0(\sigma^0_{k})\}, \quad  k\geq 0.
\end{align}
For some $k$ s.t.  $\sigma^0_k<\infty$, we have  for all  $t\in\llbracket\sigma^0_k +1,\sigma^0_{k+1}\rrbracket$:
\begin{align} \label{eq-PconstX0}
\bP^{\mbox{{\tiny{MKV}}}}(t) &= \;  \bC^{t-\sigma^0_{k}} \bP^{\mbox{{\tiny{MKV}}}}(\sigma^0_k)  
+ \sum_{s=0}^{t-\sigma^0_k-1} \bC^s \bP^0(t-1-s) \\
&= \;   \bC^{t-\sigma^0_{k}} \bP^{\mbox{{\tiny{MKV}}}}(\sigma^0_k)  + (I_K - \bC^{t-\sigma_k^0})(I_K- \bC)^{-1} \bP^0(\sigma_k^0),
\end{align}  
where $I_K$ denotes the identity matrix in $\R^{K\times K}$. 
Notice that $I_K-C$ is invertible. Indeed, the Frobenius norm of $\bC$ satisfies: $\| \bC\|$ $=$ 
$\sqrt{\sum_{k,\ell} |\mu_k \bolc_{k\ell}|^2}$ $\leq$ $\sqrt{\sum_k \mu_k \sum_{\ell} |\bolc_{k\ell}|}$ by noting that 
$0\leq |\mu_k|^2 \leq \mu_k \leq 1$, $0\leq |c_{k\ell}|^2 \leq c_{k\ell} \leq 1$. It follows that $\|\bC\|$ $<$ 
$\sum_k \mu_k$ $=$ $1$, and thus for $\|\bC^n\|$ $\leq$ $\|C\|^n$ $<$ $1$ for all $n$ $\in$ $\N$, which shows that $I_K -\bC$ is invertible. 
Hence,  at time $t$ between two changes of opinion $\sigma^0_k$ and $\sigma^0_{k+1}$, an individual uses as an approximation of  the average opinion the barycentric (in the case $K$ $=$ $1$) combination between $\bP^{\mbox{{\tiny{MKV}}}}(\sigma^0_{k})$ and $(I_K- \bC)^{-1} \bP^0(\sigma^0_{k})$, with the weight $\bC^{t-\sigma^0_k}$. As  $\bC^{t-\sigma^0_k}\underset{t\rightarrow \infty}{\longrightarrow} 0$,  it follows that the main influencers opinions  will have a greater and greater impact with time, until the next change of average opinion of main influencers occurs, i.e.,  at $\sigma^0_{k+1}<\infty$.

Also notice that on the set  $\{\sigma^0_{k+1}=\infty\}$, the major influencers average opinion remains constant infinitely after $\sigma^0_{k}$, and  the proportion of opinion $1$ in the population  will converge:
\beqs
\lim_{t\to\infty}\bP^{\mbox{{\tiny{MKV}}}}(t){\bf 1}_{\{\sigma^0_{k+1}=+\infty\}} &=& 
(I_K- \bC)^{-1} \bP^0(\sigma_k^0) {\bf 1}_{\{\sigma^0_{k+1}=+\infty\}},
\enqs 
The factor $(I_K-\bC)^{-1}$ means that the social influence will tend to a limit proportion that is an exaggerated version of the major influencers weighted distribution, which is consistent with the phenomenon of herd behavior. This is easier to see for $K=1$. In this case, $\bC$ is a scalar, and, if there was no social influence, i.e. $\bC=0$, the limit would be $\bolc^0 \boX^0(\sigma^0_{k})$, but with social influence, i.e. $\bC>0$, the more people are dragged to the limit distribution, the more people will be attracted to it.  


\begin{Example}[Analysis of Example \ref{ex:cycles}]This example has two classes and one main influencer switching opinions each $T$ units of time. With the notation of this section, we have
\[
\bC=\begin{pmatrix}
a & 0\\
b & b
\end{pmatrix}, \quad
\boC^0= \begin{pmatrix}
a\\
b
\end{pmatrix}
\]where $a=\mu_1/2$ and $b=\mu_2/3$. Therefore:
\[
\bC^t=
\begin{pmatrix}
a^t & 0\\
\sum_{k=1}^t a^kb^{t-k} & b^t
\end{pmatrix}.
\]
 Let  $\sigma=\sigma^0_k$, for some $k$, i.e., any element of the collection of switching times $\{kT,k\in\mathbb N\}$. Then, for $t<T$, the expression \eqref{eq-PconstX0} writes:
\begin{equation}
\begin{cases}
    P_1(\sigma+t) \; = \;  a^t P_1(\sigma)+\frac{1-a^t}{1-a}aX^0(\sigma)\\
    P_2(\sigma+t) \; = \;  \sum_{s=1}^t a^ s b^{t-s} P_1(\sigma)+b^t P_2(\sigma)+\frac{1-b^t}{1-b}bX^0(\sigma). 
\end{cases}    
\end{equation}
\end{Example}

Further investigations below will rely on some specific assumptions for the dynamics of the main influencer opinion process $\boX^0$: (i) in Markov assumption in Subsection \ref{subs:X0Markov} and (ii) deterministic, regular fluctuations in Subsection \ref{subs:X0det}. In the sequel, when in the single class case $K$ $=$ $1$, we omit the bold notation for an element in $\R^K$ $=$ $\R$, and write e.g. $\bP$ $=$ $P$.

\subsection{Markov assumption for the main influencer opinion dynamics}\label{subs:X0Markov}

In this subsection, we assume the linear dynamics together with the following Assumption: 

\begin{Assumption}
The opinion process of the main influencers, i.e., the process $\boX^0$ is  a Markov chain. We denote by $\mathcal I^0=\{0,1\}^{M^0}$ the state space of this process, and  $Q=(q_{i,j})_{(i,j)\in\mathcal I^0\times \mathcal I^0 }$ its transition matrix.
\end{Assumption}

\begin{Theorem}\label{thmconvergence}
We suppose that the transition matrix $Q$ is irreducible with invariant distribution $\pi$ and that $\boX^0 \sim $ Markov$ (\pi,Q)$. Consider the  process $\bP$  solution of \eqref{PMKVlinmulti} starting from $\bP(0)=0$. Then $\bP(t)$ converges in distribution to 
\begin{align}\label{eq-Pinfty}
\bP_\infty & = \;   \sum_{t=1}^\infty \bC^{t-1} \boCc^0  \boY(t), 
\end{align}
with $\boY \sim$ Markov$ (\pi,\widehat Q)$ where the transition matrix $\widehat Q$ has the  components $\widehat q_{i,j}= \frac{\pi(j)}{\pi(i)} q_{j,i}$, for $i,j\in\mathcal I^0$. Furthermore, the following hold:
\begin{itemize}
\item[(1)]  The expected value of $\bP_\infty$ coincides with the one of the weighted average opinion of the main influencers:
 \begin{align}
\E[\bP_\infty] &= \;  (I_K- \bC)^{-1}  \E[\bP^0(0)].
 \end{align}
In the single class case $K$ $=$ $1$, $\bC$ $=$ $c$ is a constant in $[0,1)$, and  the variance is explicitly expressed as: 
 \beqs
 \text{Var}(P_\infty) &=& 
 \frac{1}{1-c^2}\left[\text{Var}(P^0(0)) +\frac{2c}{1-c}\text{Cov}(P^0(0), P^0(\Gamma))\right], 
 \enqs
 where $\Gamma$ is a geometric random variable, with parameter $(1-c)$ and independent of the process $P^0$.
\item[(2)]   Suppose that $\{\boX^0(t),t\in \mathbb N\}$ are independent,  then the cumulants $\varphi_n$ of $\bP_\infty$  satisfy, in the single class case $K$ $=$ $1$, the relation 
\beqs
\varphi_n= \frac{1}{1-c^n}\varphi^0_n,
\enqs 
with $\varphi_n^0$ being the $n$-th cumulant of $\bP^0(0)$. In particular, for $ n=1,2,3$, the cumulants are the $n$-th central moments.
\end{itemize}
\end{Theorem}

The proof is available in Subsection \ref{proof:thmconvergence}.
\bigskip

\begin{Remark}
A particular case is when  there is a single main influencer with opinion following a Bernoulli(1/2) process (i.e., $X^0$ follows a sequence of i.i.d. Bernoulli random variables with parameter $1/2$), and there is a single class $K$ $=$ $1$ of individuals, assigning the weight $c^0=1-c$ to the main influencer's opinion. Then, the distribution of the limiting opinion can be linked to the Bernoulli convolutions. Indeed, a  linear transformation  of $P_\infty$, namely  $\beta:= \frac{c}{1-c}(2P_\infty -1)$ leads to $\beta=\sum_{k=1}^\infty c^k B(k)$ where $B$ is a sequence of independent random variables, whose distribution satisfies $\P[B(k) = -1]$ $=$ $\P[B(k) = 1]$ $=$ $1/2$. The Bernoulli convolution with parameter $c\in(0,1)$ is the probability measure $\nu_c$ on $\mathbb R$ that is the distribution of the random variable $\beta$.
The study of Bernoulli convolutions has a long history, going back to Jessen and Wintner \cite{JessWint35} that proved that the distribution $\nu_c$ is either singular or absolutely continuous with respect to the Lebesgue measure  (see e.g. Varj\'u \cite{Varj16} for a recent survey of the topic).   In our setting,  the following conclusions can be drawn from the existing literature: 
\begin{itemize}
\item[-] If $c=\frac{1}{2}$ then $P_\infty\sim \text{Uniform}[0,1]$.
\item[-] If $c=\frac{1}{3}$ then the limit distribution is the Cantor distribution.
\item[-] More generally, for any $c<\frac{1}{2}$ the distribution of $P_\infty$ is supported by a generalised Cantor set (hence is singular w.r.t.  the Lebesgue measure).
\item[-] Almost every $c$ in the interval $(\frac 12,1) $ gives rise to an absolutely continuous distribution such that its density is in $L^2(\mathbb R)$ (Solomyak \cite{Solo95}); examples of $c\in (\frac 12,1)$ with $\nu_c$ singular are also known (e.g. Erd\"os \cite{Erdo39}). 
\end{itemize}
These mathematical results are connected to the phenomenon of social inertia that will be analysed in more details in the next section.  When $c=\frac{1}{3}$, the social influence is weaker than when  $c=\frac{1}{2}$, and thus the major influencer has a stronger influence. Therefore, when the major influencer changes his opinion, the public opinion will evolve much faster in the direction of this new opinion, which will translate to larger jumps in the public opinion's evolution than for $c=\frac{1}{2}$. When the major influencer has opinion $1$, for instance, the next value of the public opinion cannot be below $\frac{2}{3}$, because even if it started from the left-most value $0$, it would instantaneously jump to $\frac{2}{3}$. Likewise, when the major influencer has opinion $0$, the next value of the public opinion cannot be above $\frac{1}{3}$, because even if it started from the right-most value $1$, it would instantaneously jump to $\frac{1}{3}$. These large jumps thus imply that the interval $(\frac{1}{3},\frac{2}{3})$ can never be visited by the public opinion, meaning that there are gaps in the public opinion's distribution support. If we carried the same analysis for $c=\frac{1}{2}$, this interval would be $(\frac{1}{2},\frac{1}{2})=\emptyset$, intuitively justifying why there is no gap in the public opinion's distribution support in this case.
\end{Remark}

\bigskip

\begin{Example}[One major influencer switching opinions at geometric times] \label{ex1}
We con\-sider the case of a single main influencer,  $\mathcal I^0=\{0,1\}$, with transition matrix $Q$ with elements  denoted by  $q_{1,0}=\alpha$ and $q_{0,1}=\beta$, and a single class $K$ $=$ $1$ of individuals. For simplicity, we consider $c^0=1-c$. 
We compute explicitly the variance of $P_\infty$ in this case. Elementary calculations lead to the stationary distribution 
$\pi$ satisfying:
\[
\pi(0) \; = \; \frac{\alpha}{\alpha+\beta},  \; \text{ and } \; \pi(1) \; = \; \frac{\beta}{\alpha+\beta}
\] 
and for $k\geq 1$:
 \begin{align*}
 \E[P^0(0)P^0(k)]&= (1-c)^2 \pi(1)q_{1,1}^{(k)}
 \end{align*}
 with the $k$-step transitions probability from state 1 to state 1:
  $$
  q^{(k)}_{1,1}=\pi(1)+\frac{\alpha}{\alpha+\beta}(1-\alpha-\beta)^k.
  $$ 
Hence 
\begin{align*}
 \E[P^0(0)P^0(\Gamma)]& = (1-c)^2\pi(1)\sum_{k=1}^\infty c^{k-1}(1-c)q_{1,1}^{(k)}\\
 & =  (1-c)^2\left(\pi(1)^2 +\pi(0)\pi(1)\frac{(1-\alpha-\beta)(1-c)}{1-(1-\alpha-\beta)c}\right).
 \end{align*}
It follows that $\text{Cov}(P^0(0),P^0(\Gamma))= \E[P^0(0)P^0(\Gamma)]-(1-c)^2\pi(1)^2=  (1-c)^2 \pi(0)\pi(1)\frac{(1-\alpha-\beta)(1-c)}{1-(1-\alpha-\beta)c}$, and $\text{Var}P^0(0)= (1-c)^2\text{Var}(X^0(0))= (1-c)^2\pi(0)\pi(1) $. Finally,  using  Theorem \ref{thmconvergence} we obtain:
 \begin{align*}
\text{Var}(P_\infty)&=\text{Var}( P^0_\infty ) \frac{1}{(1-c^2)}\left( \frac{1+(1-\alpha-\beta)c}{1-(1-\alpha-\beta)c}\right).
 \end{align*}
We observe that if $\alpha+\beta<1$, then the variance of $P_\infty$ is larger than the variance of the weighted main influencers opinion under the stationary distribution, $P^0_\infty$ and is increasing in the social influence factor $c$. 

Coming back to Example \ref{ex:eco-chambers}, we notice that the two classes follow trajectories that are not depending of the other class (the two communities do not interact). Hence we can apply the result of this section for each class, paying nevertheless attention to the fact that class 2 has non constant coefficients of influence, and also not summing up to 1. In this example $\alpha=\epsilon=1-\beta$. We obtain:
\begin{itemize}
    \item[-] for class 1, the limit opinion distribution has expected value $1-\epsilon$, close to 1 (as $\epsilon$ assumed small) and variance $3\epsilon(1-\epsilon)/4$, that is, close to zero. The interpretation is that in the long run, the opinion 1 will be predominant.
    \item[-] for class 2, it can be easily seen that the macro-scale function equals: $\phi_2(p,x)=p_2-\nu p_2 1_{x=1}$. Hence, the opinion process in class 2, $P_2(t)$ remains constant whenever $X^0(t-1)=0$ and decreases by a proportion $\nu$ of its level, whenever $X^0(t-1)=1$. Thus, as the process $X^0$ visits the level 1 infinitely often, the process $P_2$ converges to 0 a.s. at $t$ goes to infinity. The opinion 1 will become non existent in the long run.  
\end{itemize}
\end{Example}

\begin{Example}[Optimal opinion diffusion]

We consider again the case of a single major influencer as above, switching opinion at geometric times and, additionally, we suppose that all individuals have the opinion 0 at time 0.  We take the point of view of a decision maker that aims at spreading the opinion 1 in the population. For achieving this objective, he is able to influence the distribution of the opinion process of the main influencer, by paying a price (e.g., influencing the intensity of the advertising of a new product, by paying a price).   

Let us  assume that at time $0$, the opinion process of the major influencer is chosen in the family of 2-state Markov chains $\{X^0_\beta,\beta\in[0,1]\}$ where $X_\beta^0$ has transition probabilities: $q_{1,0}=\alpha$ and  $q_{0,1}=\beta$,  where $\alpha\in(0,1)$ is a fixed constant that is not influenced by the decision maker. Let us introduce the cost function:
\begin{align*}
 f(\beta)=\sum_{s=1}^\infty\rho^s \E[P(s)-X_\beta^0(s)+\theta \beta],
\end{align*}
where $\rho\in(0,1)$ is a fixed constant (the discounted rate), $\theta$ a nonnegative constant 
We assume that the benefit at time $t$, $P(t)$  is reduced by the cost of advertising: whenever the main influencer will diffuse opinion 1 there is a cost of  $1$. We also incorporate a cost for not advertising, as this time will be used to spreading the opposite opinion 0 instead. This is a cost of missed opportunity, assumed proportional to the average time the main influencer spends in state 0 and given by $-\theta(1-\beta)=\theta\beta-\theta$ (we exclude the constant $\theta$ from the function $f$ as it does not impact the optimisation problem considered).  
We consider the  problem:
\begin{align*}
\max_{\beta\in[0,1]}\;& f(\beta).
\end{align*}
Assuming  $X^\beta_0(0)=0$ a.s. and $\xi_n=0$  a.s. for all $n$ (so that $P(0)=0$), we obtain
\begin{align*}
\E[P(t)-X_\beta^0(t)]&=(1-c)  \sum_{s=0}^{t-1} c^s q^{(t-1-s)}_{0,1}(\beta) -q^{(t)}_{0,1}(\beta).
\end{align*}
 It is a standard result that the $t$-step transition probability $q_{0,1}^{(t)}$ has the expression:
$$
q_{0,1}^{(t)}(\beta):=\P[X^0_\beta(t)=1|X^0_\beta (0)=0] \; = \; \frac{\beta}{\alpha+\beta}\left(1-(1-\alpha-\beta)^t\right),\; t\in \N.
$$
We study the different cases:
\begin{itemize}
\item[I)]If $\alpha+\beta=1- c$, we obtain:
\begin{align*}
\E[P(t)-X_\beta^0(t)]
&=  \frac{\beta}{\alpha+\beta}\left[(1-c) \left(\frac{1-c^t}{1-c}-tc^{t-1}\right)-\left(1-c^t\right)\right]\\
&= -\beta tc^{t-1}\leq 0
\end{align*}
and we obtain:
\[
 g(\beta):=\sum_{s=1}^\infty\rho^s \E[P(s)-X^\beta_0(s)]=-\frac{\beta}{1-\rho c}\left(\frac{\rho}{1-\rho (1-\alpha-\beta)}\right)
\]
\item[II)]If $\alpha+\beta\neq 1-c$, then:\\
\begin{align*}
\E[P(t)-X_\beta^0(t)]&=  \frac{\beta}{\alpha+\beta}\left[(1-c) \left(\frac{1-c^t}{1-c}-\frac{(1-\alpha-\beta)^{t}-c^t}{(1-\alpha-\beta)-c}\right)+(1-\alpha-\beta)^t-1\right]\\
&=  \frac{\beta}{1-\alpha-\beta-c}\left( c^t-(1-\alpha-\beta)^{t}\right).
\end{align*}
Therefore:
\begin{align*}
 g(\beta):=\sum_{s=1}^\infty\rho^s \E[P(s)-X^\beta_0(s)]=& \frac{\beta}{1-\alpha-\beta-c}\left(\frac{1}{1-\rho c}-\frac{1}{1-\rho(1-\alpha-\beta)}\right)\\
&=-\frac{\beta}{1-\rho c}\left( \frac{\rho}{1-\rho (1-\alpha-\beta)}\right),
\end{align*}
that is the same expression as in the case I).
\end{itemize}
After rewriting $ g(\beta)=\frac{1}{1-\rho c}(\frac{1}{1+M\beta}-1)$, with $M=\rho/ (1-\rho+\alpha\rho)>0$
it may be checked that on [0,1], $ g\leq 0$ with $g(0)=0$, and decreasing convex. As $f(x)= g(x)+x \frac{\theta}{1-\rho}$, it follows that if $\theta$ sufficiently high,  $f$ may become nonnegative and increasing on an interval $[a,1]$, and remain negative on $[0,a)$; otherwise (i.e. $\theta$ not sufficiently high) $f$ remains negative. Hence the condition $f(1)> 0$ (that is equivalent to: $ \theta> \frac{\rho(1-\rho)}{(1-\rho c)(1+\alpha \rho)}$) guarantees that $a<1$ and hence is optimal to choose $\beta=1$, while $f(1)< 0$ means that $f(x)<0$ for all $x\in[0,1]$, hence optimal $\beta$ is 0. 

We can conclude that the advertising strategy depends on how large is the cost of missed opportunity $\theta$:
\begin{itemize}
\item[-] If $ \theta< \frac{\rho(1-\rho)}{(1-\rho c)(1+\alpha \rho)}$, then 
$$ 
\argmax_{\beta\in[0,1]} f(\beta)=0
$$ and the opinion 1 will not diffuse. Indeed, the major influencer will never adopt the opinion 1, and neither will the population of individuals. 
\item[-] If $ \theta> \frac{\rho(1-\rho)}{(1-\rho c)(1+\alpha \rho)}$ then it is optimal to use the the maximum influence $\beta =1$.
\item[-] If  $\theta= \frac{\rho(1-\rho)}{(1-\rho c)(1+\alpha \rho)}$ the decision maker is indifferent between no influence ($\beta =0 $) and maximum influence ($\beta =1$).
\end{itemize}

\end{Example}

\subsection{Regular and deterministic fluctuations for the main influencers opinion process}\label{subs:X0det}

In this subsection, we still assume a  linear macro-scale function, a single class $K$ $=$ $1$ of individuals, together with the following  assumption: 
  
\begin{Assumption} \label{Hypregfluc}
There is a single major influencer with opinion process $\boX^0$ defined by:
\beqs 
\boX^0(t)=
\begin{cases}
1 \quad   {\rm  for }  \quad  t\in \llbracket 2kT, (2k+1)T\llbracket, \; k\in \N, \\
0 \quad   {\rm  otherwise,}
\end{cases}
\enqs 
for some given $T$ $\in$ $\N^*$. 
\end{Assumption}

Assumption \ref{Hypregfluc} means that  the major influencer starts with opinion $1$ and alternates between states  $0$ and $1$; it spends  $T$ consecutive units times in a state, before moving to the other state. Equivalently, $\boX^0$ is characterized by the relation $\sigma^0_k$ $=$  $kT$ for $k\in\N$, where $\{\sigma^0_k,k\in\N\}$ is the sequence of stopping times defined in \eqref{eq-sigma0}.  We have the following quantitative behavior about the mean-field approximation  sequence  
$\{P^{\mbox{{\tiny{MKV}}}}(t), t \in \N\}$ of the population's average opinion, starting from $P^{\mbox{{\tiny{MKV}}}}(0)$ $=$ $0$.


\begin{Theorem} \label{theofluctu} 
The sequence $\{P^{\mbox{{\tiny{MKV}}}}(t), t \in \N\}$  is  increasing in the intervals $\llbracket 2kT, (2k+1)T\llbracket$, and decreasing in the intervals $\llbracket (2k+1)T, (2k+2)T\llbracket$, and:
\[
P^{\min}_{\infty}:=\liminf_{t\to\infty} P^{\mbox{{\tiny{MKV}}}}(t)= \frac{c^T}{1+c^T}\frac{\bc^0}{1-c}
\]

\[
P^{\max}_{\infty}:=\limsup_{t\to\infty} P^{\mbox{{\tiny{MKV}}}}(t)= \frac{1}{1+c^T}\frac{\bc^0}{1-c}.
\] 
\end{Theorem}


The proof is found in Subsection \ref{proof:theofluctu}


\vspace{3mm}

\begin{Remark}We observe that the higher the cycle length $T$ is, less frequent oscillations are, but with a higher amplitude $P^{max}_\infty-P^{min}_\infty$ (as $P^{min}_\infty= c^T P^{max}_\infty$).  
Another key role is played by the social influence parameter $c$. In Subsection \ref{subs:NumX0det}, we show that increasing $c$ produces opposite effects on the fluctuations of the public opinion, that we call snowball effect and social inertia.
\end{Remark}

We  apply this result to a practical problem, where a trade off between frequency and amplitude of oscillations is searched.

\begin{Example}[Optimal cycle of switching]
Assume that an agent is selling products subject to trends (e.g. clothes, decoration, furniture, etc). Let us consider  two styles of product, style $A$ and style $B$, that are sold in alternation (i.e., same duration, non-overlapping  periods). For an individual, being in state $x\in\{A,B\}$ is interpreted as  using a product of style $x$. The process $X^0(t)$ keeps track of the current style (the only that  can be bought at time $t$). An individual in state $x$ at time $t$ was satisfied with the product style at time $t-1$; whenever he/she become insatisfied, she/he will buy the current product and start using it at time $t+1$, hence adopt the state $X^0(t)$ at time $t+1$. 

The agent collaborates with the major influencer and chooses the period $T$ optimally, to incite people to buy products as often as possible. Let us assume that the gain of the agent at time $t\in\N^\star$ is 
$Q \vert P^{\mbox{{\tiny{MKV}}}}(t)-P^{\mbox{{\tiny{MKV}}}}(t-1)\vert$ for some $Q\in\R_+$. It simply means that his gain is proportional (with factor $Q$) to the proportion $\vert P^{\mbox{{\tiny{MKV}}}}(t)-P^{\mbox{{\tiny{MKV}}}}(t-1)\vert$ of people who bought a product (the current type) at time $t-1$.  The agent's gain is measured by his asymptotic average gain per time, that is:
\beqs 
V(T) &:=& \lim_{t\rightarrow \infty}\frac{1}{t}\sum_{s=1}^t Q\vert P^{\mbox{{\tiny{MKV}}}}(s)-P^{\mbox{{\tiny{MKV}}}}(s-1)\vert,
\enqs 
where we stress the dependence on the duration $T$ of the change cycle. 

From Theorem \ref{theofluctu}, it is clear, given the oscillating nature of $P^{\mbox{{\tiny{MKV}}}}$ every $T$ period, that  for any $m$ $\in$ $\N^\star$, 
\beqs 
&&\sum_{s=1}^{2mT} \vert P^{\mbox{{\tiny{MKV}}}}(s)-P^{\mbox{{\tiny{MKV}}}}(s-1)\vert\\
&=& \sum_{k=0}^{m-1} \Big( \big(P^{\mbox{{\tiny{MKV}}}}((2k+1)T)-P^{\mbox{{\tiny{MKV}}}}(2kT) \big) \big)+ 
\big( P^{\mbox{{\tiny{MKV}}}}((2k+1)T)-P^{\mbox{{\tiny{MKV}}}}(2(k+1)T)\Big)\\
&=& \sum_{k=0}^{m-1} \Big( P^{max}_k-P^{min}_k+  P^{max}_{k}-P^{min}_{k+1}\Big),
\enqs 
using the same notation as in \eqref{eq:Pmin} and \eqref{eq:Pmax}. Therefore, 
\beqs 
V(T)&=&\frac{1}{2T}\lim_{m\rightarrow \infty}\frac{1}{m}\sum_{k=0}^{m-1} \Big( P^{max}_k-P^{min}_k+  P^{max}_{k}-P^{min}_{k+1}\Big) \; = \; \frac{1}{2T}2 (P^{max}_\infty-P^{min}_\infty)\\
&=&\frac{1}{T}\frac{1-c^T}{1+c^T}\frac{\bc^0}{1-c}.
\enqs 
The goal of the agent is to maximize $V(T)$ over the duration $T$. The derivative of $T$ $\mapsto$ $V(T)$ is given by:
\beqs 
V'(T)= -\frac{1}{T^2}\frac{1-c^T}{1+c^T}+\frac{1}{T}\frac{-2\ln(c)c^T}{(1+c^T)^2}. 
\enqs
The sign of $V'(T)$ clearly is the same as the sign of $c^{-T}(1+c^T)^2 T^2 V'(T)$, which writes as
\beqs 
-(c^{-T}-c^{T})-2\ln(c)T, 
\enqs
and it is  negative if and only if
\beqs 
c^{-T}-c^{T}+2\ln(c)T\geq 0.
\enqs
By setting $y=-\ln(c)T$, this is written as 
\beqs 
e^{y}-e^{-y}-2y\geq 0, 
\enqs 
which is equivalent to $sh(y)\geq y$, which is true for all $y\geq 0$.  
Therefore, $V'(T)\leq 0$, and thus $V$ is maximized by $T=1$ (the smallest possible value of $T$). This means that the optimal behavior for the major influencer is to change his mind as often as possible. In this case, the population will asymptotically oscillate between the proportions $c\frac{\bc^0}{1-c^2}$ and $\frac{\bc^0}{1-c^2}$, and the optimal average gain per time is
\beqs 
V^\star=V(1)=\frac{\bc^0}{1+c}.
\enqs 
An interesting qualitative consequence of this formula is that the optimal gain of the agent is proportional to $\bc^0$, the influence factor of the major influencer, and inversely proportional to $1+c$, in particular it is decreasing in the population's influence factor. An intuition for these properties is that it is the major influencer's influence which leads the sudden trend change, while the population's influence acts as an impediment to the change, working against the major influencer's influence. Indeed, for each individual there is  a trade-off between following the major influencer and staying with the majority of the population, and this trade-off will lean toward staying with the population if its influence factor $c$ is large, and toward following the major influencer if his influence factor $\bc^0$ is large.

\end{Example}

\section{Numerical illustrations} \label{sec:num}

In this section, we graphically illustrate via some examples the phenomena theoretically derived in previous sections, namely on the one hand, propagation of chaos and on the other hand, snowball effect and social inertia.

For our study, we  shall assume that all thresholds $(\mrs_n)$ are uniform on $[0,1]$. There is a single class of agents ($K=1$), and 
the major influencer has deterministic opinions, either fixed, or following periodic cycles. Consequently, the mean field approximation will also be deterministic in all the examples.  For given parametrizations of the model, all figures from this section follow the convention that the mean-field public opinion process is displayed in green, while the samples of the 
$N$-agent public opinion process $(\bP^N(t))_t$ are displayed in red. The different instances correspond to different public opinions trajectories, as the model’s features vectors $(\chi_n)$ are re-sampled from one instance to another one. 


\subsection{Some phenomena related to propagation of chaos}

We  first study a toy  example designed to illustrate the various phenomena that are specifically related to propagation of chaos.

We assume that $\xi_n\sim \text{Bernoulli}(1/2)$ for all $n\in \llbracket 1,N\rrbracket$,  there is a single major influencer always in state $1$, and  
$\bolc^0=(1-c)/2$, where $c$ is any positive real. One can easily show that the corresponding mean-field dynamic is constantly equal to $\frac{1}{2}$, thus its graphical representation is the horizontal line with equation $y=\frac{1}{2}$. Hence, this example is convenient to study the propagation of chaos,  by simply comparing the $N$-agent dynamics to this horizontal line.

We will illustrate four phenomena relative to the propagation of chaos:
\begin{itemize}
    \item {\it The impact of the population size $N$:} we compare the complete information to the mean-field dynamics for different sizes of the population $N$, illustrating that the larger the population, the better the mean-field approximation.
    \item {\it The role of the social influence on the divergence regimes of the $N$-agent and the mean-field dynamics:}  the level of social influence factor $c$ determines if the macro-scale function $\boldsymbol\phi$ is contracting or not, creating convergence versus divergence of opinions towards the mean-field dynamics.  
     as $t\rightarrow \infty$.
    \item {\it The impact of common surveys on propagation of chaos:} we illustrate that in the case of  common surveys,  the propagation of chaos is actually driven by $M$, the size of the population sample. The larger $M$ is, the better the mean-field approximation will be.
    \item {\it The impact of independent surveys on propagation of chaos:} in the case of  independent surveys instead of common surveys, the propagation of chaos is actually not driven by the sample size, and the mean-field approximation is still good even if $M=1$. 
\end{itemize}

\begin{figure}[H]
\centering
\includegraphics[height=4.5cm,width=.32\textwidth]{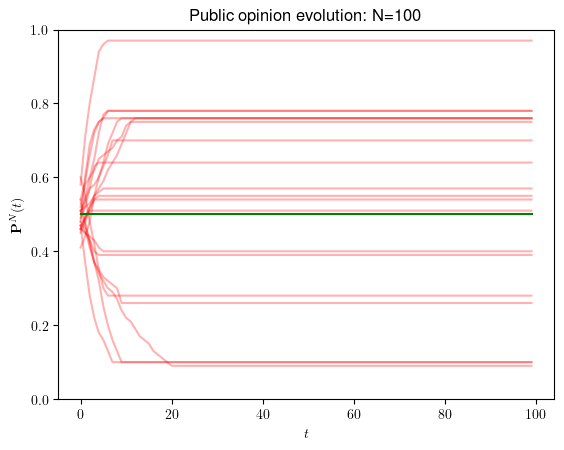}\hfill
\includegraphics[height=4.5cm,width=.32\textwidth]{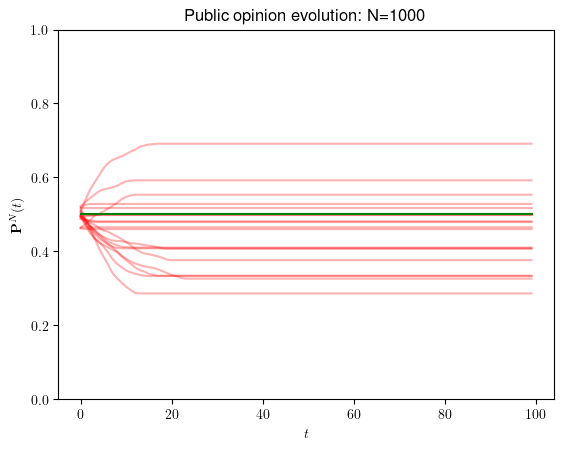}\hfill
\includegraphics[height=4.5cm,width=.32\textwidth]{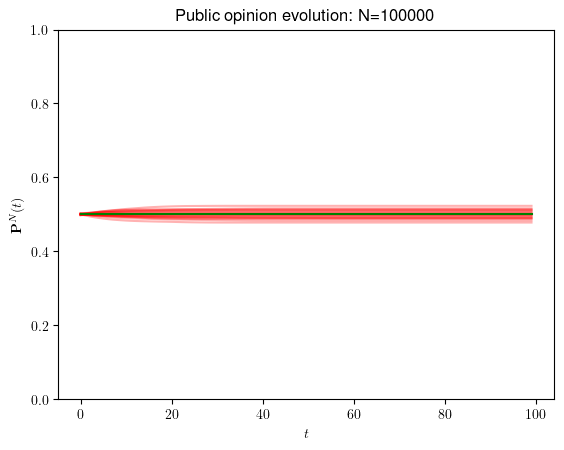}
\caption{Propagation of chaos. Left: $N=100$. Middle: $N=1000$, Right: $N=1000000$}
\label{fig:prop}
\end{figure}

In Figure \ref{fig:prop}, we plotted the $N$-agent dynamics with complete information for different values of $N$ ($N=100$, $1000$, $100000$), with the fixed parameters $c=0.9$,  
and we observe that the ``band'' formed by the red curves is thinner and thinner around the green line as $N$ is larger, which means that the mean-field approximation is better for large $N$.

Let us now fix $N=100000$. In principle, that should lead  to a reasonably good approximation. However, that approximation is generally guaranteed to be good only at the beginning (i.e. $t=0$), and may or may not diverge as $t$ increases, depending on the constant $c$. (recall that the contracting constant is $K_\phi$ $=$ $c$ in this linear case example).

We now vary $c$ by taking the values $c=0.9$, $1$, and $1.1$. 
In Figure \ref{fig:div}, we see that 

\begin{itemize}
    \item for $c<1$, the ``red band'' stays close to the green line. 
\item for $c=1$, the ``red band'' linearly diverges far away from green line. 
\item for $c>1$, the ``red band'' exponentially diverges far away from green line, which confirms the relevance of $e^{(K_{\phi}-1)T}$ term in the bound of Theorem \ref{theo-glob}.
\end{itemize}  
These experiments confirm the bound in  Theorem \ref{theo-glob}. 

\begin{figure}[H]
\centering
\includegraphics[height=4.5cm,width=.32\textwidth]{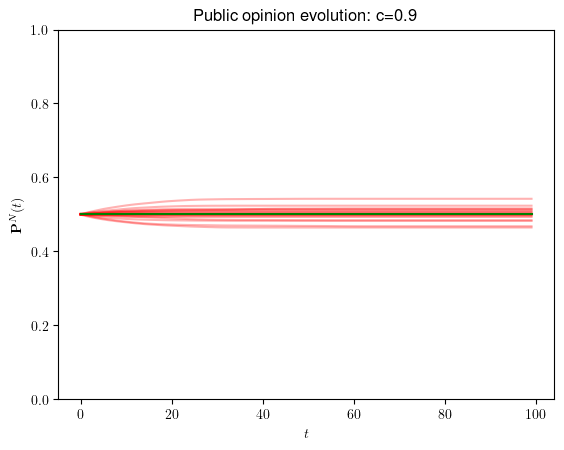}\hfill
\includegraphics[height=4.5cm,width=.32\textwidth]{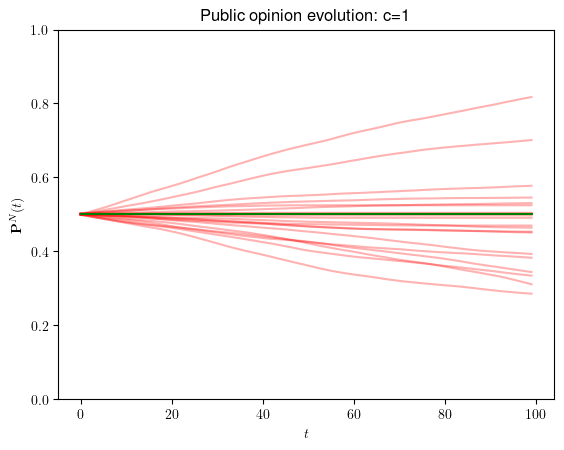}\hfill
\includegraphics[height=4.5cm,width=.32\textwidth]{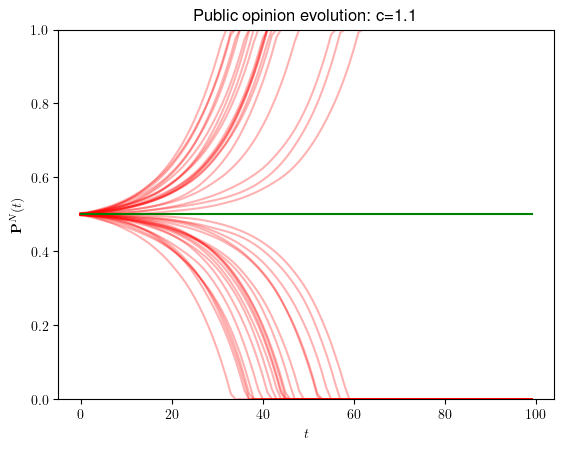}
\caption{Divergence regimes. Left: $c=0.9$. Middle: $c=1$, Right: $c=1.1$} 
\label{fig:div} 
\label{fig:div}
\end{figure}

We also note the ``smoothness''  of the $N$-agent dynamics. Indeed, in this complete information setting and with a fix opinion of the main influencer, the randomness only affects the features vectors, that are sampled randomly at the initial time. Conditionally on these, the curves of the opinions dynamics are deterministic. 
The situation becomes very different if we introduce incomplete information, such as the use common or independent surveys.

\begin{figure}[H]
\centering
\includegraphics[height=4.5cm,width=.32\textwidth]{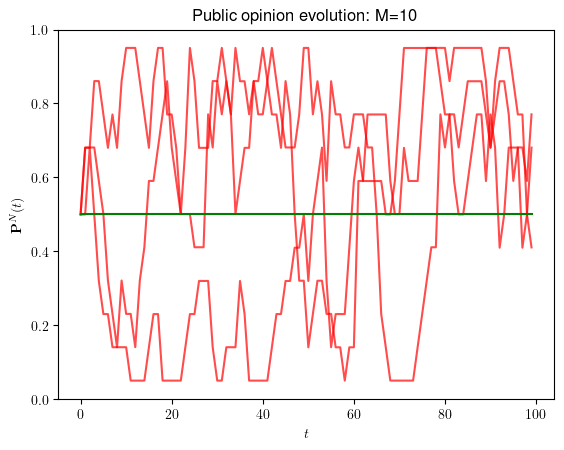}\hfill
\includegraphics[height=4.5cm,width=.32\textwidth]{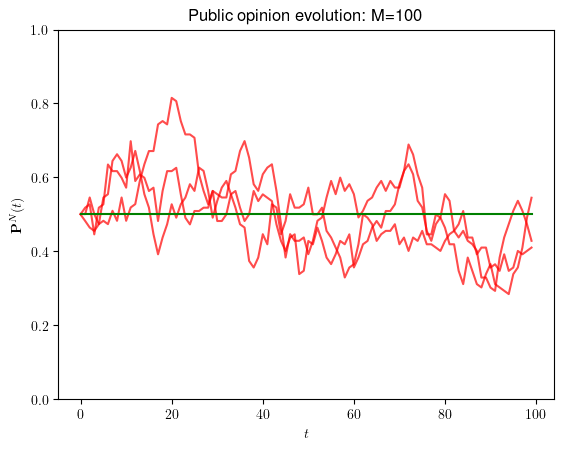}\hfill
\includegraphics[height=4.5cm,width=.32\textwidth]{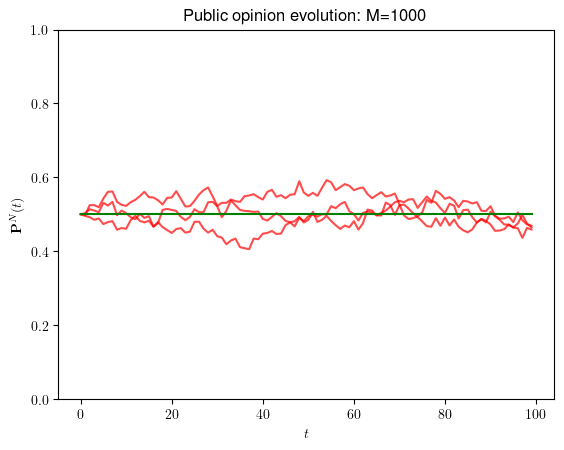}
\caption{Common survey impact. Left: $M=10$. Middle: $M=100$, Right: $M=1000$}
\label{fig:common}
\end{figure}

Figure \ref{fig:common} illustrates the $N$-agent dynamic with common surveys, comparatively to the mean-field. We used the fixed parameters $N=1000000$, and  the parameter $M$ (size of the common sample) took valued $M=10$, $100$, and $1000$.  We observe that indeed, using surveys creates randomness in the time evolution of the $N$-agent dynamic, as randomness does not only occur for sampling the features of the individuals but also at each time $t$, when the surveyed groups are re-sampled.
Another phenomena that we clearly observe is that the quality of the mean-field approximation depends upon $M$, the size of the surveyed group: the larger it is, the better the mean-field approximation, which is consistent with Theorem 
\ref{theo-glob}.

\begin{figure}[H]
\centering
\includegraphics[height=4.5cm,width=.32\textwidth]{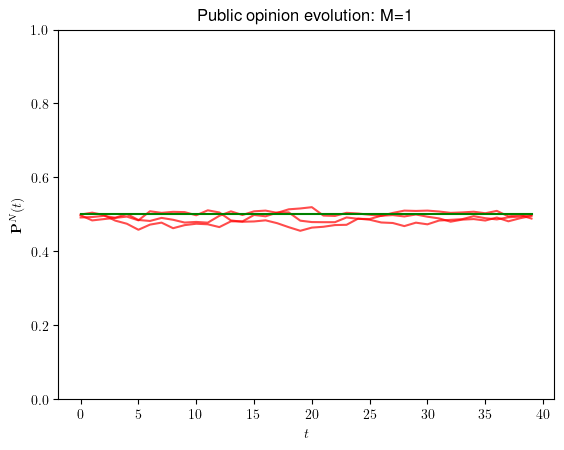}\hfill
\includegraphics[height=4.5cm,width=.32\textwidth]{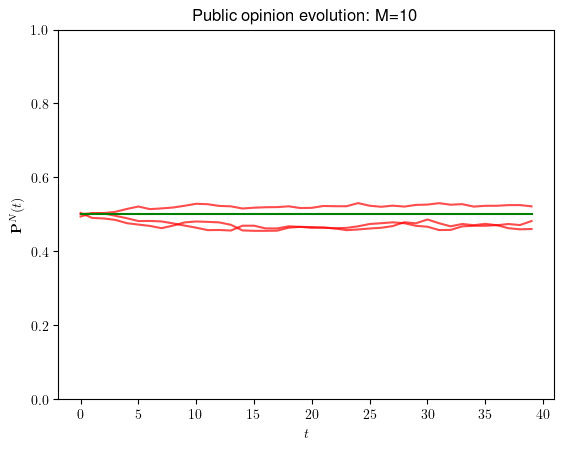}\hfill
\includegraphics[height=4.5cm,width=.32\textwidth]{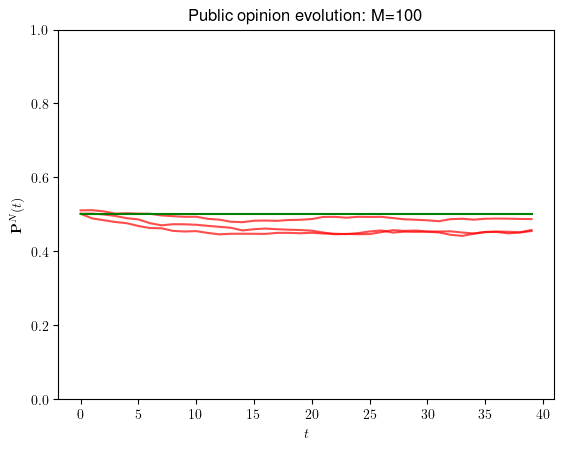}
\caption{Independent surveys impact. Left: $M=1$. Middle: $M=10$, Right: $M=100$}
\label{fig:indep}
\end{figure}

Figure \ref{fig:indep} illustrates the $N$-agent with independent surveys, comparatively to the mean-field public opinion. We used the fixed parameters $N=10000$, $c=0.9$, while the parameter $M$ is varying with the values $M=1$, $10$, $100$. We observe that although there is still noise in the $N$-agent's dynamic (again due to the fact that surveyed groups are re-sampled at each time $t$), the situation qualitatively differs from the dynamic with common survey in the sense that the quality of the mean-field approximation here seems not be affected a lot  by the size of $M$: even for $M=1$, the mean-field approximation is good. Although this is not stated in Theorem \ref{theo-glob}, this phenomenon can be intuitively explained by the fact  that even though each individual samples only one opinion when $M=1$, overall, $N$ individuals have been surveyed, and even though each individual is only influenced by the opinion of one of these $N$ individuals, the public opinion is influenced by all of them. By contrast, in the case of common surveys with $M=1$, at each time the whole population is influenced by the opinion of a single individual. This is why in the independent survey case, even with $M=1$, a form of law of large numbers will take place and make the dynamic close to its mean-field approximation. 

\begin{figure}[H]
\centering
\includegraphics[height=5cm,width=.32\textwidth]{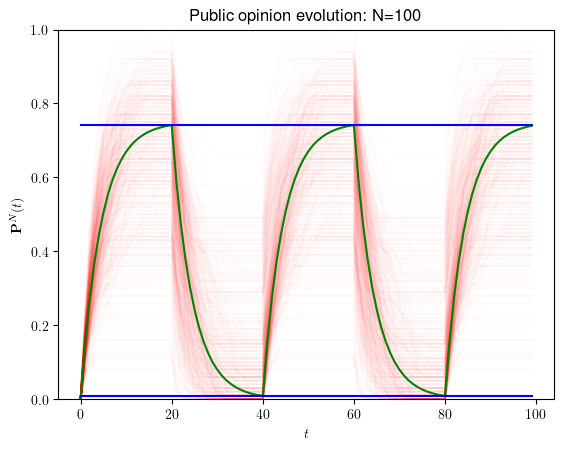}\hfill
\includegraphics[height=5cm,width=.32\textwidth]{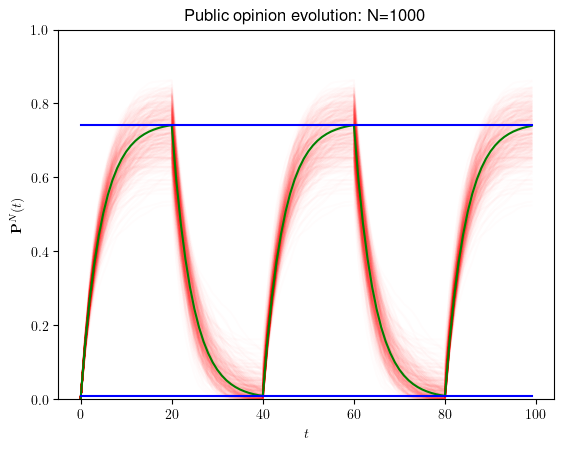}\hfill
\includegraphics[height=5cm,width=.32\textwidth]{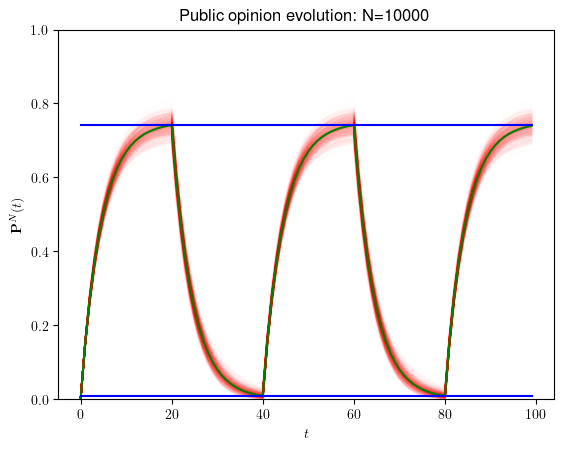}
\caption{Propagation of chaos. Left: $N=100$. Middle: $N=1000$, Right: $N=10000$}
\label{fig:fluct-prop}
\end{figure}

Next, we illustrate the propagation of chaos for the problem in Section \ref{subs:X0det}. 
In Figure \ref{fig:fluct-prop}, we plot the $N$-agent and mean-field public opinion's dynamics with fixed parameters $c=0.8$, $\bolc^0=0.15$, $T=20$ and the parameter $N$ taking the values $N=100$, $1000$, $10000$. The blue levels represent the maximum and minimum proportion of opinion $1$ as given by $P^{\min}_\infty$  and $P^{\max}_\infty$ in Theorem \ref{theofluctu}. One can clearly see that the larger $N$ is, the thinner is the ``red curved band'' (representing the approximate proportion of public opinion) around the green curve levels. This means that the larger $N$ is, the closer the $N$-agent dynamic with complete information will be to the mean-field dynamic, and therefore the better the mean-field approximation will be.

\subsection{Graphical illustration of the snowball effect and social inertia
}
\label{subs:NumX0det}

We consider  the problem from Section \ref{subs:X0det}, and illustrate the two following phenomena:
\begin{itemize}
        \item {\it Snowball effect:} The fact that as the social influence factor $c$ gets larger, the maximal proportion of opinions $1$ increases.
    \item {\it Social inertia:} The fact that the social influence factor $c$ acts as a brake on the evolution of public opinion, leading to two phenomena: by increasing $c$, (1) the transition from one minimum to one maximum (and vice versa) will be slower and slower, and (2) the local maximum and minimum of the oscillations of the public opinion will converge more and more slowly to their limit values.
\end{itemize}

Let us fix $N$ $=$ $10000$ to have a reasonably good approximation and let us focus on illustrating the properties of the mean-field dynamic in this example.

\begin{figure}[H]
\centering
\includegraphics[height=5cm,width=.32\textwidth]{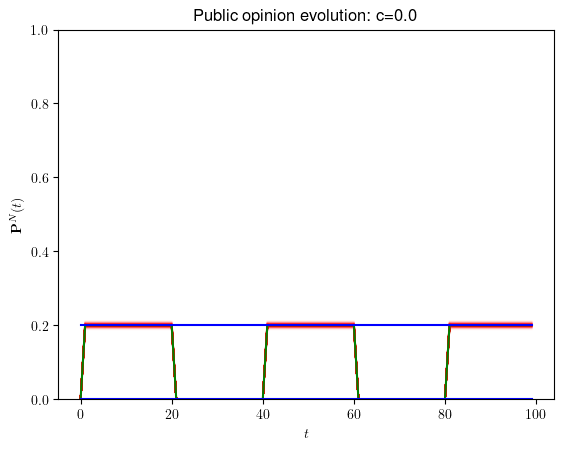}\hfill
\includegraphics[height=5cm,width=.32\textwidth]{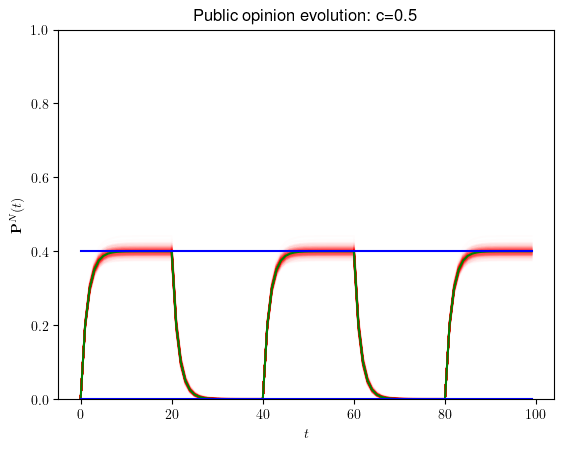}\hfill
\includegraphics[height=5cm,width=.32\textwidth]{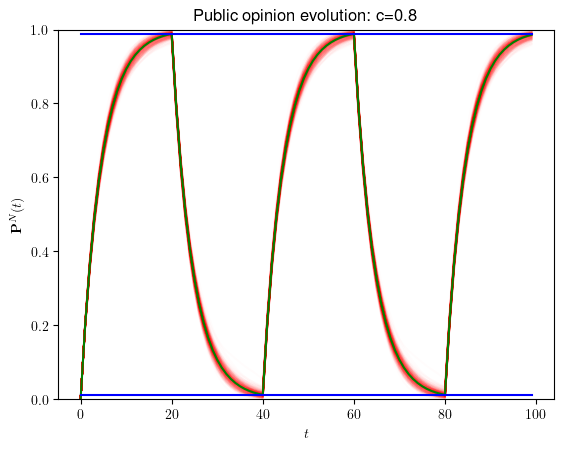}
\caption{Snowball effect. Left: $c=0$. Middle: $c=0.5$, Right: $c=0.8$}
\label{fig:snowball}
\end{figure}

One of the phenomena we mentioned in Section \ref{subs:X0det} was the ``snowball effect'': when the major influencer has opinion $1$, the proportion of opinion $1$ in the public opinion is more dragged toward $1$ if there is social influence ($c>0$) compared to a population only influenced by the major influencer ($c=0$), and it is more and more dragged toward $1$ as $c$ is large. 
In Figure \ref{fig:snowball}, we plotted the $N$-agent and mean-field dynamics for the fixed parameters $N=10000$, $\bolc^0=0.15$, $T=20$, and the parameter $c$ taking the values $c=0$, $0.5$, $0.8$. One can see that the local maximum of the oscillations (corresponding to how close to $1$ the public opinion is dragged when the major influencer has opinion $1$) are indeed higher for larger $c$.

We also illustrate another impact of the social influence parameter $c$, namely social inertia. Notice that as $c$ is increased, not only the maximal public opinion increases (snowball effect), but also, the {\it speed} at which it converges to its maximum, and then back to its minimum, seems to {\it decrease}. To better illustrate this phenomenon, we performed the same experiments, still with $N=10000$ and $T=20$, with $c$ taking the values $c=0.3$, $0.5$, $0.8$, also adjusting $\bolc^0$ to maintain the same $0.9$ distance between the limit maximum and minimum of the oscillations.

\begin{figure}[H]
\centering
\includegraphics[height=5cm,width=.32\textwidth]{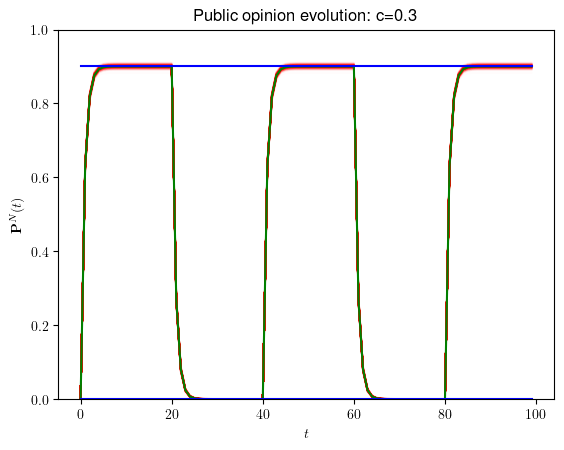}\hfill
\includegraphics[height=5cm,width=.32\textwidth]{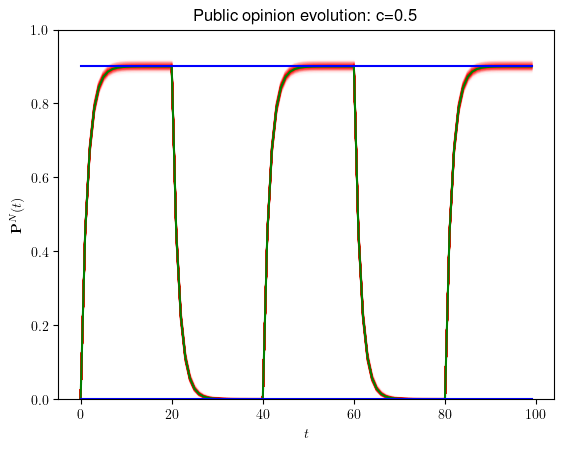}\hfill
\includegraphics[height=5cm,width=.32\textwidth]{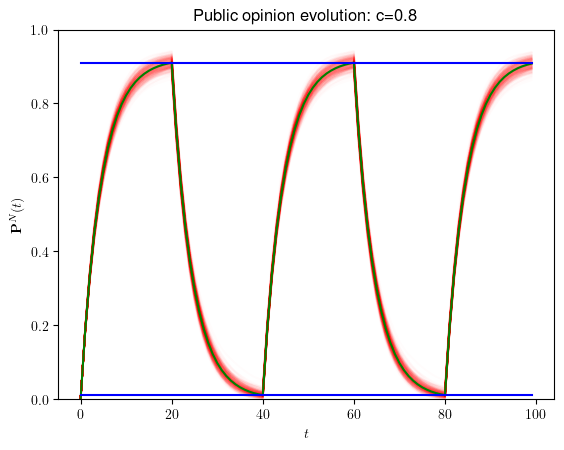}
\caption{Social inertia: speed of convergence to local maximum/minimum. Left: $c=0.3$. Middle: $c=0.5$, Right: $c=0.8$}
\label{fig:inertia-speed}

\end{figure}

In Figure \ref{fig:inertia-speed}, we see that indeed, $c$ also affects the ``shape'' of the oscillations: for $c$ small, the oscillations have a ``square'' shape, and for $c$ larger, they get more and more ``triangular'', meaning that for $c$ small the public opinion almost immediately switches from its minimum to its maximum (and vice versa), while for $c$ large it slowly and smoothly transitions between these values, as if there was an opposite force retaining it and slowing its transitions.

\begin{figure}[H]
\centering
\includegraphics[height=5cm,width=.32\textwidth]{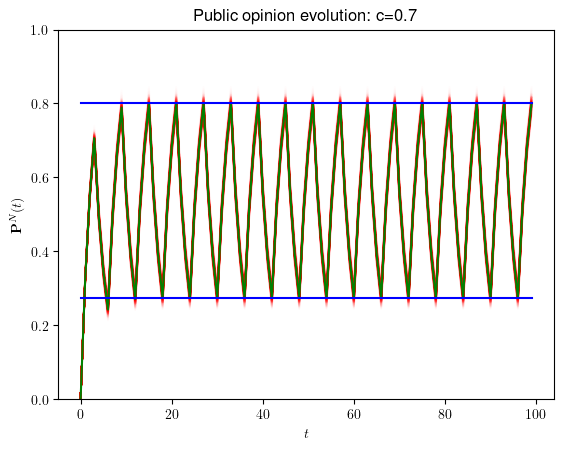}\hfill
\includegraphics[height=5cm,width=.32\textwidth]{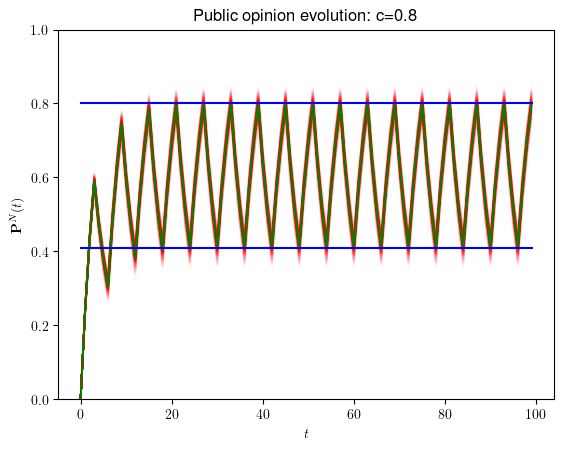}\hfill
\includegraphics[height=5cm,width=.32\textwidth]{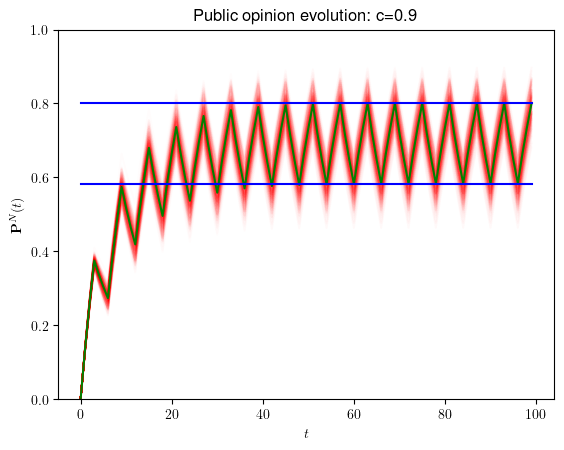}
\caption{Social inertia: convergence to limit oscillations. Left: $c=0.7$. Middle: $c=0.8$, Right: $c=0.9$}
\label{fig:inertia-limit}

\end{figure}

This suggests a last phenomenon that we illustrate here: if $c$ gets even larger, the transition from minimum to maximum might be so slow that the public opinion doesn't have the time to get to its asymptotic maximum from the first oscillation, and might require several oscillations before converging to its limit oscillations.
In Figure \ref{fig:inertia-limit} we plotted the $N$-agent and mean-field public opinion dynamics for fixed parameters $N=10000$, and $c$ taking the values $c=0.7$, $0.8$, $0.9$, and adjusting $\bolc^0$ to have a constant asymptotic upper bound. We see that indeed, for larger $c$, as $c$ increases, the number of initial oscillations required to converge to the  asymptotic oscillations (between the two horizontal blue lines) is greater and greater.
Consequently, the social influence factor $c$ has two effects that are, interestingly, ``opposite'' to each other in a sense: it increases the amplitude of the oscillations (snowball effect), while reducing the speed of the public opinion's evolution (social inertia).

\appendix

\section{Proofs of the main results}\label{sec-Proofs}

\subsection{Proof of Theorem \ref{theo-loc}}\label{proof:theo-loc}
\noindent{\bf Case of McKean-Vlasov approximation.} 
Let us first introduce some notation:
$\mathcal X_t:=\mathcal X\times(\{0,1\}\times\{0,1\}^{M_0})^t\text{ for }t\geq 1$, where we recall  that $\mathcal X=\{0,1\}\times\Kc\times \R^{M_0}\times\R$. For $n\in\llbracket 1,N\rrbracket $, and for $t\geq 0$, denote 
\begin{equation}\label{eq-H}
H_n(t):=(\chi_n, (B_n(s),\boX^0(s))_{s\leq t})
\end{equation}
 i.e., the vector containing the features of individual $n$ and her/his information up to time $t$.  We set $H_n(-1):=\chi_n$.

Let us prove by induction that for all $t\in\N$, there exists a function $f_t:\mathcal X_t\rightarrow \{0,1\}$ such that
\begin{align}\label{eq-fn}
\tilde X^{}_n(t)=f_t(H_n(t-1)),\quad a.s.
\end{align}
Let us fix $(x, k,c^0,\mrs)\in\mathcal S$. For $t=0$, it is clear that $f_0((x, k,c^0,\mrs))=x$ satisfies the wished property. Assuming that for some $t\in\N$, such function $f_t$ exists, let us build a suitable function $f_{t+1}$. For  $h \in \mathcal X_t$ let  $\widehat{h}$ be vector containing the first  elements of $h$, such  that $h$ has the representation $(\widehat{h},(b(t),x^0(t)))$, where $(b(t),x^0(t))\in \{0,1\}\times\{0,1\}^{M_0}$.  We assume $(x, k,c^0,\mrs)\in \mathcal S$ are the first four elements in  $\widehat h$. Also, we denote by $g_t(h)= \sum_\ell c_{k,\ell} \P[f_t(H_n(t-1))=1,\kappa_n=\ell]+c^0 \cdot x^0(t).$ 
By assumption,  \eqref{eq-fn} holds, and hence,  
\[
g_t(H_n(t) )=\sum_{\ell\in\mathcal K}c_{\kappa_n,\ell}\P[\tilde X^{}_n(t)=1,\kappa_n=\ell\mid \Fc^0_t]+\bc^0 \cdot \boX^0(t)
\]  
so that $\{g_t(H_n(t) )>\mrs_n\}$ is the event that individual $n$ updates her/his opinion to 1 at time $t+1$, assuming  $t$ is an influence time with corresponding  available information $H_n(t)$. Now, we define
\beqs 
f_{t+1}(h(t))&:=&(1-b(t))f_t(\widehat{h})
 + b(t) {\bf 1}_{g_t(h)> \mrs}. 
\enqs 
Hence:
\begin{align*}
f_{t+1}(H(t))=&(1-B_n(t))f_t(H(t-1))+ \;  B_n(t) {\bf 1}_{g_t(H_n(t)) > \mrs_n}\\
=&(1-B_n(t)) \tilde X^{}_n(t)+ B_n(t) {\bf 1}_{\bc_{\kappa_n}\cdot(\P(\tilde X^{}_n(t),\kappa_n=\ell\mid \Fc^0_t))_{\ell\in\Kc}+\bc^0 \cdot \boX^0_t> s_n} \\
 = &  \tilde X^{}_n(t+1), 
\end{align*}
which concludes the induction. 
We then have:
\begin{align}
&\E\;\Big[\big\vert \E\left [\tilde X^{}_n(t)) {\bf 1}_{\kappa_n=k}\mid \Fc^0(t)\right]-\frac{1}{N} \sum_{i\leq N} 
\tilde X^{}_i(t){\bf 1}_{\kappa_i=k}\big\vert^2\Big] \\
=& \; \E\;\Big[ \vert \E\left [f_t(H_n(t-1) ){\bf 1}_{\kappa_n=k}\mid \Fc^0(t)\right ]-\frac{1}{N}\sum_{i\leq N} f_t(H(t-1)){\bf 1}_{\kappa_i=k}\vert^2\Big] \\
=& \;  \E\;\Big[\big \vert \E\left[f_t(\chi_n,(B_n(s), x^0(s))_{s\leq t-1}))\right ]-\frac{1}{N}\sum_{i\leq N} f_t(\chi_i,(B_i(s), x^0(s))_{s\leq t-1})\big\vert^2\Big\vert_{(x^0(s):=\boX^0(s)),\forall s\leq t-1} \Big]\\
=& \; \E\;\Big[\text{Var}\Big(\frac{1}{N}\sum_{i\leq N} \tilde X^{}_i(t){\bf 1}_{\kappa_i=k}\;\big|\;\Fc^0(t) \Big) \Big], 
\end{align}
and consequently:
\beqs 
\E\Big[ \sum_{k\in\Kc} \big \vert \tilde{P}_{n,k}(t)&-&\frac{1}{N}\sum_{i\leq N} \tilde X^{}_i(t) 1_{\kappa_i =k} \big\vert^2 \Big]\\
&=&\sum_{k\in\Kc} \E\Big[\big\vert \E[\tilde X^{}_n(t)) {\bf 1}_{\kappa_n=k}\mid \Fc^0(t)]-\frac{1}{N}\sum_{i\leq N} 
\tilde X^{}_i(t){\bf 1}_{\kappa_i=k}\big\vert^2 \Big]\\
&=&\sum_{k\in\Kc}  \E\;\Big[\text{Var}\Big(\frac{1}{N}\sum_{i\leq N} \tilde X^{}_i(t){\bf 1}_{\kappa_i=k}\;\big|\;\Fc^0(t) \Big) \Big]\\
&=&\E\Big[\text{Var}\Big(\frac{1}{N}\sum_{i\leq N}\tilde X^{}_i(t)\;\big|\;\Fc^0(t) \Big) \Big].
\enqs
It follows that 
\beqs
\E\Big[ \sum_{k\in\Kc}   \vert P_{k}^{MKV}(t)-\frac{1}{N}\sum_{i\leq N} \tilde X^{}_i(t) 1_{\kappa_i=k} \vert^2 \Big] &\leq &
\frac{1}{N}\E\Big[\text{Var}\Big( f_t(\chi_i,(B^i_s, x^0_s)_{s\leq t-1})\Big)_{x^0_s:=\boX^0_s,s\leq t-1} \Big] \\
&\leq& \frac{1}{4N}. 
\enqs 
by the property of the Bernoulli distribution whose variance is bounded by $1/4$, and thus by Cauchy-Schwarz inequality: 
\begin{align} \label{inegloc}
\Ic_{\ell oc}^N(t) & \leq \;   \frac{\sqrt{K}}{2} \frac{1}{\sqrt{N}}. 
\end{align}

\vspace{5mm}

\noindent{\bf Case of common or independent sample groups approximation:} 
The cases of common or independent sample groups approximation are much simpler and are both treated in the exact same way. Let us thus only prove the result for the case of independent sample groups. We have 
\beqs
\E\Big[ \sum_{k\in\Kc} \big|  \tilde P_{n,k}(t) &-&\frac{1}{N}\sum_{i\leq N} \tilde X^{}_i(t) {\bf 1}_{\kappa_i=k} \big|^2 \Big]\\
&=& \E\Big[ \sum_{k\in\Kc}  \big| \frac{1}{M}\sum_{i\in I^n_t} \tilde X^{}_i(t)-\frac{1}{N}\sum_{i\leq N} 
\tilde X^{}_i(t)\big|^2  {\bf 1}_{\kappa_i=k}  \Big]\\
&=&\frac{1}{M^2}\E\Big[ \sum_{k\in\Kc}  \E\Big[\big| \sum_{i\in I^n_t} x_i-\frac{M}{N}\sum_{i\leq N} x_i\big|^2 \Big]_{x_i=\tilde X^{}_i(t),i\leq N}{\bf 1}_{\kappa_i=k}  \Big]\\
&=&\frac{1}{M^2}\E\Big[\text{Var}(Z_{p,N,M})_{p:=\frac{1}{N}\sum_{i\leq N} \tilde X^{}_i(t)}\Big]
\enqs 
where, for $p\in [0,1]$, $Z_{p,N,M}$ is a hypergeometric random variable with parameters $(p,N,M)$ consisting in counting the number of individuals in state $1$ when drawing $M$ of them without replacement in a population of $N$ individuals containing a proportion $p$ of individuals in state $1$. By the classical properties of the hypergeometric distribution, we then deduce that
\beqs 
\E\Big[ \sum_{k\in\Kc} \big|  \tilde P_{n,k}(t) -\frac{1}{N}\sum_{i\leq N} \tilde X^{}_i(t) {\bf 1}_{\kappa_i=k} \big|^2 \Big] &=& 
\frac{1}{M^2}\E\Big[\Big(p(1-p) \frac{M(N-M)}{N-1}\Big)_{p:=\frac{1}{N}\sum_{i\leq N} \tilde X^{}_i(t)}\Big] \\
&\leq &\frac{1}{4M}\frac{N-M}{N-1}, 
\enqs 
which gives the required result by Cauchy-Schwarz inequality. 

\subsection{Proof of Theorem \ref{theo-glob}}\label{proof:theo-glob}

Notice that the macro-scale error writes:
\beqs 
\Ic^N(\tilde{\bP},T) &=&\sum_{k\in\Kc} \E\Big[\frac{1}{N}\sum_{n= 1}^N \big\vert \tilde X^{}_n(T) - X^{N}_n(T)\big\vert{\bf 1}_{\kappa_n=k} \Big]\\
&=&\frac{1}{N}\sum_{n= 1}^N \sum_{k\in\Kc} \E\Big[ \big\vert \tilde X^{\tilde \bP}_n(T) - X^{N}_n(T)\big\vert{\bf 1}_{\kappa_n=k} \Big].
\enqs 

\noindent{\bf Case of the McKean-Vlasov approximation.} 
In the case of the McKean-Vlasov approximation, all individuals $n\in \llbracket 1,N\rrbracket$ use the same approximations to form their opinion, given by the vector, $\tilde\bP$ $=$ $\bP^{\mbox{\tiny{MKV}}}=(P^{\mbox{\tiny{MKV}}}_1,\dots. P^{\mbox{\tiny{MKV}}}_K)$.

Again from the i.i.d. assumption of the features vectors  $(\chi_n)_n$,  $\E[|X^N_n(t)-\tilde X^{}_n(t)\vert {\bf 1}_{\kappa_n=k}]$ does not depend on $n$, and hence we can alternatively write the macro-scale error function as:
\begin{align} \label{eq-INdelta}
\Ic^N({\bP}^{MKV},T) 
&= \;  \sum_{k\in\Kc} \Delta_k(t),
\end{align} 
where we have denoted $\Delta_k(t):=\E[|X^N_n(t)-\tilde X^{}_n(t)\vert {\bf 1}_{\kappa_n=k}]$ for $k\in\Kc$.  

From \eqref{eq-general-dyn} and \eqref{dynXtildep}, we have 
\beqs 
|X^{N}_n(t+1)- \tilde X^{}_n(t+1)| = (1-B_n(t))|X^{N}_n(t)-\tilde X^{}_n(t)| +B_n(t) |{\bf 1}_{S_n^N(t)>\mrs_n}- {\bf 1}_{\tilde S_n(t)>\mrs_n}|.  
\enqs 

It follows that for any $n$ $\in$ $\lb 1,N\rb$, and $k$ $\in$ $\Kc$,
\begin{align} 
\Delta_k(t+1)&= \; (1-h)\Delta_k(t)+h \E[|{\bf 1}_{\bc_{\kappa_n}\cdot \bP^N(t) + \bc^0_n\cdot \boX^{0}(t)>\mrs_n}- 
{\bf 1}_{ \bc_{\kappa_n} \cdot  \bP^{\mbox{\tiny{MKV}}}(t) + \bc^0_n\cdot \boX^{0}(t)>\mrs_n}|{\bf 1}_{\kappa_n=k}]
\nonumber \\
&= \; (1-h)\Delta_k(t)+h \E\left [\left({\bf 1}_{M_k^N(t)+\bc^0_n\cdot \boX^{0}(t) >\mrs_n} - {\bf 1}_{m^N_k(t) +\bc^0_n\cdot \boX^{0}(t)>\mrs_n}\right){\bf 1}_{\kappa_n=k}\right]  \label{eq-deltak}
\end{align}  
where we have used the notation: $m^N_k(t)=\min[\bc_k\cdot \bP^N(t); \bc_k\cdot  \bP^{\mbox{\tiny{MKV}}}(t)]$, and $M^N_k(t)=\max[\bc_k\cdot \bP^N(t); \bc_k\cdot  \bP^{\mbox{\tiny{MKV}}}(t)]$. 
Further, we notice that 
\beqs
\frac{1}{N}\sum_{i=1}^N {\bf 1}_{M^N_k(t)+ \bc^0_i\cdot \boX^{0}(t)  >  \mrs_i}{\bf 1}_{\kappa_i=k} &=&  \max\big[ \Phi^N_k(\bP^N(t), \boX^{0}(t))\;,\;\Phi^N_k(\bP^{\mbox{\tiny{MKV}}} (t), \boX^{0}(t)) \big]
\enqs
and 
\beqs
\frac{1}{N}\sum_{i=1}^N {\bf 1}_{m^N_k(t)+ \bc^0_i\cdot \boX^{0}(t)> \mrs_i}{\bf 1}_{\kappa_i=k} &=&  \min\big[ \Phi^N_k(\bP^N(t), \boX^{0}(t))\;,\; \Phi^N_k(\bP^{\mbox{\tiny{MKV}}} (t), \boX^{0}(t)) \big],
\enqs
as the empirical average functions $\Phi^N_k$ are nondecreasing. Using these expressions and the symmetry between the individuals in the population, we obtain:
\beqs 
&&\sum_{k\in\Kc}\E\left [\left({\bf 1}_{M_k^N(t)+\bc^0_n\cdot \boX^{0}(t) >\mrs_n} - {\bf 1}_{m^N_k(t) +\bc^0_n\cdot \boX^{0}(t)>\mrs_n}\right){\bf 1}_{\kappa_n=k}\right]\\
&=&\sum_{k\in\Kc}\frac{1}{N}\sum_{i=1}^N\E\left [\left({\bf 1}_{M_k^N(t)+\bc^0_i\cdot \boX^{0}(t) >\mrs_i} - {\bf 1}_{m^N_k(t) +\bc^0_i\cdot \boX^{0}(t)>\mrs_i}\right){\bf 1}_{\kappa_n=k}\right]\\
&=&\E\left[|  \Phi^N( \bP^N (t), \boX^{0}(t))-\Phi^N(\bP^{\mbox{\tiny{MKV}}} (t), \boX^{0}(t))|\right]\\
&\leq&\E\; |  \Phi^N( \bP^N (t), \boX^{0}(t))-\boldsymbol \phi(\bP^N(t), \boX^{0}(t))|+     \E\;|  \Phi^N( \bP^{\mbox{\tiny{MKV}}}(t), \boX^{0}(t))- \boldsymbol \phi(\bP^{\mbox{\tiny{MKV}}} (t), \boX^{0}(t))|\\
&& + \E\;| \boldsymbol  \phi ( \bP^N (t), \boX^{0}(t))-\boldsymbol \phi (\bP^{\mbox{\tiny{MKV}}} (t), \boX^{0}(t))|\\
&\leq& 2 \E\; ||  \Phi^N-\boldsymbol \phi ||+ \E\;| \boldsymbol  \phi ( \bP^N (t), \boX^{0}(t))-\boldsymbol \phi (\bP^{\mbox{\tiny{MKV}}} (t), \boX^{0}(t))|\\
&\leq &  2 \E\; ||  \Phi^N-\boldsymbol \phi ||+ K_\phi \E\;|\bP^N (t)- \bP^{\mbox{\tiny{MKV}}} (t))|\\
&\leq &  2 \E\; ||  \Phi^N-\boldsymbol \phi ||+ K_\phi \E\;\sum_{k\in\Kc} \Big(\Big |P_k^N (t)- \frac{1}{N}\sum_{i=1}^N 
\tilde X_i^{} (t) {\bf 1}_{\kappa_i=k}\Big |+\Big |\frac{1}{N}\sum_{i=1}^N \tilde X_i^{} (t) {\bf 1}_{\kappa_i=k}- \bP^{\mbox{\tiny{MKV}}} (t))\Big |\Big)\\
&\leq& 2 \E\; ||  \Phi^N-\boldsymbol \phi ||+ K_\phi \Big(\Ic^N(\bP^{\mbox{\tiny{MKV}}},t)+\Ic_{\ell oc}^N(\bP^{\mbox{\tiny{MKV}}},t)\Big). 
\enqs 
Using the inequality above and \eqref{eq-INdelta}-\eqref{eq-deltak},  we obtain:
\beqs 
 \Ic^N(\bP^{\mbox{\tiny{MKV}}},t+1)&=&\sum_{k\in\Kc} \Delta_k(t+1)  \\
& \leq& (1-h)\Ic^N(\bP^{\mbox{\tiny{MKV}}},t)+h \left(K_\phi \;\Ic^N(\bP^{\mbox{\tiny{MKV}}},t) 
+  K_\phi \Ic_{\ell oc}^N(\bP^{\mbox{\tiny{MKV}}},t)  + 2 \E\; ||  \Phi^N-\boldsymbol \phi ||\right)\\
 &\leq &(1-h+hK_\phi )\Ic^N(\bP^{\mbox{\tiny{MKV}}},t)+ h \Big(\frac{K_\phi \sqrt{K}}{2\sqrt{N}}+ 2\E\; ||  \Phi^N-\boldsymbol \phi|| \Big), 
 \enqs
hence by induction:
 \beqs 
 \Ic^N(\bP^{\mbox{\tiny{MKV}}},t)&\leq &\frac{1-(1-h+hK_\phi )^{t}}{1-K_\phi} 
\Big(\frac{K_\phi \sqrt{K}}{2\sqrt{N}} + 2\E\; ||  \Phi^N-\boldsymbol \phi|| \Big) \\
& \leq &  \frac{1 - K_\phi^t}{1-K_\phi}  
\Big(\frac{K_\phi \sqrt{K}}{2\sqrt{N}} + 2\E\; ||  \Phi^N-\boldsymbol \phi|| \Big), 
 \enqs
 and the result follows.
 \ep

\vspace{5mm}

 \noindent{\bf Case of common or independent sample groups approximation:} 
 
 We adapt the previous proof and treat these two cases together. At the difference with the McKean-Vlasov approximation,  here individuals do not use the same public opinion approximation for updating opinions: a priori, at any point in time $t$, we have $\tilde{\bP}_n(t)\neq \tilde{\bP}_m(t)$ for $n\neq m$. 
 
 However, in the case of common (respectively   independent) $M$-sample approximations, at any time $t\in \mathbb N^*$, the random variables $I(t)$ (respectively $(I_n(t))_{n\in\llbracket 1,N\rrbracket }$) are independent form $(H_n(t))_{n\in\llbracket 1,N\rrbracket }$ (see \eqref{eq-H} for the definition of these random variables).  Furthermore, for $t\neq s$, we have that $I(t)$ and $I(s)$ are independent   (respectively $(I_n(t))_{n\in\llbracket 1,N\rrbracket }$ and $(I_n(s))_{n\in\llbracket 1,N\rrbracket }$ are independent). That means that the surveyed groups at time $t$,  used for updating the opinions at time $t$, are drawn independently of the history of all  the observations in the population. Moreover, from one individual to the other, either the surveyed groups are the same, or they are independent (depending on the framework considered). Consequently, in this case as well, $\E[|X^N_n(t)-
 \tilde X^{}_n(t)\vert {\bf 1}_{\kappa_n=k}]$ does not depend on $n$ and we can follow  the steps of the previous proof with some small adaptation.
 
 We use the notation: $m^N_{n,k}(t)=\min[\bc_k\cdot \bP^N(t); \bc_k\cdot  \tilde \bP_n(t)]$, and $M^N_{n,k}(t)=\max[\bc_k\cdot \bP^N(t); \bc_k\cdot  \tilde \bP_n(t)]$, as these quantities are not the same for all individuals $n$. 
 However, by the independence properties mentioned above, we get for $n\neq r$:
 \beqs 
R_k(t):&=&\P\left (\mrs_n- \bc^0_n\cdot \boX^{0}(t) \in (m^N_{n,k}(t)  , M^N_{n,k}(t)] ; \kappa_n=k\right)\\
&=& \P\left (\mrs_r- \bc^0_r\cdot \boX^{0}(t)  \in (m^N_{r,k}(t) , M^N_{r,k}(t)] ; \kappa_r=k\right)
\enqs
(hence this quantity does not depend on the individual). Therefore, similarly to \eqref{eq-deltak}, we obtain 
\beqs 
\Delta_k(t+1)
&=&(1-h)\Delta_k(t)+h \E[|{\bf 1}_{\bc_{\kappa_n}\cdot \bP^N(t) + \bc^0_n\cdot \boX^{0}(t)>\mrs_n}- 
{\bf 1}_{ \bc_{\kappa_n} \cdot \tilde  \bP_n(t) + \bc^0_n\cdot \boX^{0}(t)>\mrs_n}|{\bf 1}_{\kappa_n=k}]\\
&=&(1-h)\Delta_k(t)+h R_k(t)
\enqs 
Again following the steeps from the previous proof, we obtain:
\beqs 
\sum_{k\in\Kc}R_k(t)&=&\sum_{k\in\Kc}\E\left [\left({\bf 1}_{M_{n,k}^N(t)+\bc^0_n\cdot \boX^{0}(t) >\mrs_n} - {\bf 1}_{m^N_{n,k}(t) +\bc^0_n\cdot \boX^{0}(t)>\mrs_n}\right){\bf 1}_{\kappa_n=k}\right]\\
&\leq& 2 \E\; ||  \Phi^N-\boldsymbol \phi ||+ K_\phi \Big(\Ic^N(\tilde \bP,t)+\Ic_{\ell oc}^N(\tilde \bP,t)\Big). 
\enqs 

Using the inequality above, \eqref{eq-deltak}, and Theorem \ref{theo-loc}, we obtain:
\beqs 
 \Ic^N(\tilde \bP,t+1)&=&\sum_{k\in\Kc} \Delta_k(t+1)=  (1-h)\Ic^N(\tilde \bP,t) +h  \sum_{k\in\Kc}R_k(t) \\
 &\leq &(1-h(1-K_\phi) )\Ic^N(\tilde \bP,t)+ h \Big(\frac{K_\phi\sqrt{K}}{2\sqrt{M}}+ 2\E\; ||  \Phi^N-\boldsymbol \phi|| \Big), 
 \enqs
hence:
 \beqs 
 \Ic^N(\tilde \bP,t+1)&\leq &\frac{1-(1-h(1-K_\phi) )^{t}}{1-K_\phi} \Big(\frac{K_\phi\sqrt{K}}{2\sqrt{M}}+ 2\E\; ||  \Phi^N-\boldsymbol \phi|| \Big), 
 \enqs
 and the result follows.

\subsection{Proof of Theorem \ref{thmconvergence}}\label{proof:thmconvergence}

We define the time reversed process $\boY_\theta$ $:=$ $(\boX^0(\theta-t))_{t\in[0,1,...,\theta]}$ with $\theta\in \mathbb N$ being a fixed time.  From standard theory on Markov processes, $\boY_\theta$ is also a Markov chain with initial distribution $\pi$ and transition matrix $\widehat Q$ (hence independent of $\theta$). Furthermore, $\pi$ is an invariant distribution for $\boY_\theta$.  It follows from \eqref{PMKVlinmulti} that for all $t$,  
$\bP(t)$ $\stackrel{law}{=}$ $\sum_{s=1}^t \bC^{s-1} \boCc^0\boY(s)$, 
where $Y$ is as in the statement, and the result in \eqref{eq-Pinfty} follows.

We now prove  Assertion (1). Let us introduce the functions $\Psi_\bz(\bolx):=\E \left[e^{\bolx \cdot \bP_\infty}|\boY(1)=\bz\right]$ and $\Psi^0_\bz(\bolx):=\E \left[e^{\bolx \cdot \bP^0(1)}|\boX^0(0)=\bz\right]$, for $\bolx\in\R^K$ and $\bz\in\mathcal I^0$.  As $\bP^0(1)$ and $\bP_\infty$ take values in $\Sc_K$, these functions take finite values. The moment generating function of $\bP_\infty$ (recall that $\boY$ follows the stationary distribution $\pi$ under $\P$) writes as
 \begin{align}\nonumber
 M(\bolx) \; := \;  \E \left[e^{\bolx \cdot \bP_\infty}\right] & = \;  \sum_{\by\in \mathcal I^0} \Psi_\by(\bolx)\pi(\by)\\ 
  &= \sum_{\by\in \mathcal I^0} \E \left[e^{  \bolx \cdot  \sum_{k=1}^\infty \bC^{k-1} \bC^0 \boY(k )}|\boY(1)=\by\right]\pi(\by)\\\label{eq-MPsiz}
 &= \sum_{\by\in \mathcal I^0}e^{   \bolx \cdot \boCc^0 \by} \E \left[e^{ \bolx \cdot \bC\bP_\infty }|\boY(0)=\by\right]\pi(\by)
  \end{align}
 and 
 \begin{align}
  \E \left[e^{ \bolx \cdot \bC \bP_\infty}|\boY(0)=\by\right]&=\sum_{\bz\in \mathcal I^0} 
  \E \left[e^{(\bC\trans\bolx)\cdot \bP_\infty}|\bY(1)=\bz,\boY(0)=\by\right]\frac{\P[\boY(1)=\bz,\boY(0)=\by]}{\P[\boY(0)=\by]} \\\label{eq-Psiz}
   &=\sum_{\bz\in \mathcal I^0} \Psi_\bz( \bC\trans\bolx)\widehat q_{\by,\bz}
  \end{align}
  so that, plugging \eqref{eq-Psiz} into \eqref{eq-MPsiz}, we obtain:
  \begin{align*}
 M(\bolx)
 &= \;  \sum_{\by\in \mathcal I^0}e^{  \bolx \cdot \boCc^0 \by} \sum_{\bz\in \mathcal I^0} 
 \Psi_\bz(\bC\trans\bolx) \left(\frac{\pi(\bz)}{\pi(\by)}q_{\bz,\by}\right)\pi(\by)\\
 &= \; \sum_{\bz\in \mathcal I^0} \pi(\bz) \Psi_\bz(\bC\trans\bolx)\sum_{\by\in \mathcal I^0}e^{  \bolx \cdot \bC^0 \by} q_{\bz,\by}\\
 &= \; \sum_{\bz\in \mathcal I^0}\pi(\bz) \Psi_\bz(\bC\trans\bolx)\Psi_\bz^0(\bolx).
 \end{align*}
We differentiate the left and right side of the above with respect to $\bolx$ and then take $\bolx=0$. This yields
\beqs
\E[\bP_\infty]=\nabla M(0) &=& \sum_{\bz\in \mathcal I^0}\pi(\bz)\left( \bC\nabla\Psi_\bz(0)+\nabla\Psi^0_\bz(0)\right) \; = \; \bC \E[\bP_\infty]+\E[\bP^0(1)],
\enqs
and thus: $\E[\bP_\infty]$ $=$ $(I_K-\bC)^{-1}\E[\bP^0(1)]$ $=$ $(I_K-\bC)^{-1}\E[\bP^0(0)]$.

\vspace{2mm}

Differentiating once more, we obtain
 \begin{align}
\E[\bP_\infty \bP_\infty\trans] \; = \; D^2 M(0) & = \; \sum_{\bz\in \mathcal I^0}\pi(\bz)\Big( \bC D^2 \Psi_\bz (0) \bC\trans  + \bC \nabla\Psi_\bz(0)\nabla\Psi^0_\bz(0)\trans \nonumber \\
& \quad \quad \quad \quad \quad + \;  \nabla\Psi^0_\bz(0)(\bC\nabla\Psi_\bz(0))\trans  \; + \; D^2\Psi^0_\bz(0)  \Big),
\end{align}
and we note that 
\beqs
\sum_{\bz\in \mathcal I^0}\pi(\bz)D^2\Psi_\bz (0) & = &  \E [\bP_\infty \bP_\infty\trans], \\ 
\sum_{\bz\in \mathcal I^0}\pi(\bz)D^2\Psi^0_\bz (0) & = &  \E [\bP^0(1)\bP^0(1)\trans] \; = \; 
\E [\bP^0(0)\bP^0(0)\trans],
\enqs  
so that 
\begin{align} \label{eq-P2infty}
 \E[\bP_\infty \bP_\infty\trans] \; = \; \bC  \E[\bP_\infty \bP_\infty\trans] \bC\trans 
 \; + \; \E [\bP^0(0)\bP^0(0)\trans] + \bC\boldsymbol{\Delta} + \boldsymbol{\Delta}\trans\bC\trans,   
\end{align}
with 
\beqs
\boldsymbol{\Delta} & := & \sum_{\bz\in \mathcal I^0}\pi(\bz)\nabla \Psi_\bz (0) \nabla \Psi^0_\bz(0).
\enqs
Now,  
\begin{align*}
\boldsymbol{\Delta}
&= \; \sum_{\bz\in \mathcal I^0}\pi(\bz)\E[\bP_\infty|\boY(1)=\bz]\E [\bP^0(1)\trans|\boX^0(0)=\bz]\\
&= \; \E\left [\left(   \sum_{k=1}^\infty \bC^{k-1} \sum_{\bx\in\mathcal I^0} 
\widehat q^{(k-1)}_{Z,\bx}\boCc^0\bolx \right)\left(   \sum_{\by\in\mathcal I^0} q_{Z,\by}\; \boCc^0 \by \right)\trans \right],
\end{align*}
where $Z$ is assumed to be an $\mathcal I^0$ valued random variable with probability mass function $\pi$ and $(\widehat q^{(k)}_{i,j})$ are the $k$-step transition probabilities associated with the matrix $\widehat Q$ and $\widehat q^{(0)}_{\bz,\bx}=1_{\bx=\bz}$, and thus 
\begin{align*}
\boldsymbol{\Delta}
&= \; \sum_{k=1}^\infty \bC^{k-1} \sum_{\bx\in\mathcal I^0}\sum_{\by\in\mathcal I^0} (\boCc^0 \bolx)\trans
\boCc^0 \by \E[ \widehat q^{(k-1)}_{Z,\bx}\times  q_{Z,\by}].  
\end{align*}
It is not hard to prove by induction that,  for $j=1,\dots,\ell$:
\[
\sum_{\bolx\in\mathcal I^0}\sum_{\by\in\mathcal I^0} \boCc^0 \bolx (\boCc^0 \by)\trans \E[ \widehat q^{(\ell)}_{Z,\bolx}\times  q_{Z,\by}] \; = \; 
\sum_{\bx\in\mathcal I^0}\sum_{\by\in\mathcal I^0}  \boCc^0 \bolx (\boCc^0 \by)\trans \E[ \widehat q^{(\ell-j)}_{Z,\bolx}\times  q^{(1+j)}_{Z,\by}] ,
\]
 in particular, taking $j=\ell=k-1$ :
\begin{align*}
\sum_{\bolx\in\mathcal I^0}\sum_{\by\in\mathcal I^0} \boCc^0 \bolx (\boCc^0 \by)\trans
\E[ \widehat q^{(k-1)}_{Z,\bx}\times  q_{Z,\by}] &=\sum_{\bolx\in\mathcal I^0}\sum_{\by\in\mathcal I^0} 
\boCc^0 \bolx  (\boCc^0 \by)\trans    q^{(k)}_{\bolx,\by}\pi(\bolx)\\
&= \; \E\left [ \boCc^0 \boX^0(0)   (\boCc^0 \boX^0(k))\trans \right]
\; = \; \E\left [\bP^0(0)\bP^0(k)\trans \right]. 
\end{align*}
 Hence:
 \begin{align} \label{Delta}
\boldsymbol{\Delta}
&= \; \sum_{k=1}^\infty\bC^{k-1} \E\big[\bP^0(0)\bP^0(k)\trans \big].
\end{align}
 
Let us now focus on the single class case $K$ $=$ $1$, and set $c$ $=$ $\bC$ $\in$ $\R$. By considering the random variable $\Gamma$ with geometric distribution of parameter $1-c$, independent of $P^0$, the expression of 
$\Delta$ in \eqref{Delta} is written as
\beqs
\Delta &=& \frac{1}{1-c} \E[P^0(0)P^0(\Gamma)]. 
\enqs
Plugging this expression  into \eqref{eq-P2infty}, we obtain
\begin{align*}
\E[(P_\infty)^2]&=\frac{1}{1-c^2} \left( \frac{2c}{1-c}\E[P^0(0)P^0(\Gamma)]+ \E [P^0(0)^2]\right). 
\end{align*}
Since $\big(\E[P_\infty]\big)^2$ $=$ $\frac{1}{(1-c)^2}\big(\E[P^0(0)]\big)^2$ $=$ 
$\frac{1}{(1-c)^2}\E[P^0(0)]\E[P^0(\Gamma)]$, one can easily verify that the variance equals indeed the stated expression, which proves assertion (1).

\vspace{2mm}

Let us now prove Assertion (2) in the single class case $K$ $=$ $1$.  
When  $\{X^0(t),t\in \mathbb N\}$ are independent, with distribution $\pi$, we have that $q_{ij}=\widehat q_{ij}=\pi(j)$. Hence, 
the expression \eqref{eq-P2infty} writes: 
\[
M(x)=\E \left[e^{cx P_\infty}\right]\E\left[e^{x(c^0X^0(1))} \right]= M(cx)M^0(x)
\]
where $M^0$ is the moment generating function of $P^0(0)=\bolc^0 \boX^0$ (note that in the case of one class, $K=1$ we have that $\mu_1=1$ and hence $\bolc^0=\boCc^0$).
The cumulants $\varphi_n$ pf $P_\infty$ are obtained from a power series expansion of the cumulant generating function $\log M(x)$. As $\log M(x)= \log M(cx)+\log M^0(x)$, it follows:
\[
 \sum_{n=1}^\infty \varphi_n\frac{t^n}{n!}=\sum_{n=1}^\infty \varphi_n\frac{(ct)^n}{n!}+\sum_{n=1}^\infty \varphi^0_n\frac{t^n}{n!},
\]
leading to $(1-c^n) \varphi_n=\varphi^0_n$.

\subsection{Proof of Theorem \ref{theofluctu}} \label{proof:theofluctu}

By  \eqref{eq-attraction},  starting from $P^{\mbox{{\tiny{MKV}}}}(0)$ $=$ $0$, we clearly see that (recall that $\bc^0$ $\geq$ $0$): 
\begin{align}\label{eq-framing-PMKV}
0 \; \leq \;  P^{\mbox{{\tiny{MKV}}}}(t) \; \leq \;  \frac{\bc^0}{1-c}, \quad \forall t\in\N.
\end{align}
Together with \eqref{PMKVlinmulti}, 
this implies   that $t$ $\mapsto$  $P^{\mbox{{\tiny{MKV}}}}(t)$ is increasing on the intervals $\llbracket 2kT, (2k+1)T\llbracket$, $k\in\N$, and decreasing on the intervals $\llbracket (2k-1)T, 2kT\llbracket$, $k\in\N^\star$. Indeed, for 
$t\in\llbracket 2kT, (2k+1)T\llbracket$, $k\in\N$, we have 
\begin{align} \label{Pcroit} 
P^{\mbox{{\tiny{MKV}}}}(t+1) &= \;  cP^{\mbox{{\tiny{MKV}}}}(t)+  \bc^0 \; \geq \;  c P^{\mbox{{\tiny{MKV}}}}(t)+(1-c)P^{\mbox{{\tiny{MKV}}}}(t) \; = \; P^{\mbox{{\tiny{MKV}}}}(t),
\end{align} 
by \eqref{eq-framing-PMKV}, while for $t$ $\in$ $\llbracket (2k-1)T, 2kT\llbracket$, we have (recall that $c$ $\leq$ $1$): 
\begin{align} \label{Pdecroit} 
P^{\mbox{{\tiny{MKV}}}}(t+1) &= \;  c P^{\mbox{{\tiny{MKV}}}}(t) \; \leq \;  P^{\mbox{{\tiny{MKV}}}}(t).
\end{align} 
The sequence 
$P^{\mbox{{\tiny{MKV}}}}$ is thus oscillating, first increasing for $T$ times, and then decreasing every $T$ times, in phase with the times when the major influencer switches opinions. The sequence 
\begin{equation}\label{eq:Pmin}
    P^{min}=(P^{min}_k)_{k\in\N}:=(P^{\mbox{{\tiny{MKV}}}}(\sigma^0_{2k}))_{k\in\N}=(P^{\mbox{{\tiny{MKV}}}}(2kT))_{k\in\N}
\end{equation}
corresponds to the sequence of successive local minima of $P^{\mbox{{\tiny{MKV}}}}$, occurring every $2T$ times, starting from $t=0$. The sequence 
\begin{equation}\label{eq:Pmax}
    P^{max}=(P^{max}_k)_{k\in\N}:=(P^{\mbox{{\tiny{MKV}}}}(\sigma^0_{2k+1}))_{k\in\N}=(P^{\mbox{{\tiny{MKV}}}}((2k+1)T))_{k\in\N}
\end{equation}
corresponds to the sequence of successive local maxima of $P^{\mbox{{\tiny{MKV}}}}$, occurring every $2T$ times, starting from $t=T$.  Moreover, by solving the recurrence relations appearing in \eqref{Pcroit}-\eqref{Pdecroit},  we  easily see   that 
\begin{equation} \label{Pkminmax}  
\begin{cases} 
P^{min}_{k+1} \; = \; c^TP^{max}_{k},  \quad k \in \N, \\
P^{max}_{k} \; = \; c^TP^{min}_{k}+(1-c^T)\frac{\bc^0}{1-c}.
\end{cases}
\end{equation}
Let us now check  by induction on $k$ that the sequences $(P_k^{min})_k$ and $(P_k^{max})_k$ are increasing: as $P^{min}_0=0$ we clearly have $P^{min}_1\geq P^{min}_0$, and thus 
\beqs 
P^{max}_{1}&=&c^TP^{min}_{1}+(1-c^T)\frac{\bc^0}{1-c}\\
&\geq& c^TP^{min}_{0}+(1-c^T)\frac{\bc^0}{1-c} \; = \; P^{max}_{0}.
\enqs 
Assume that for $k\in\N$ we have $P^{min}_{k+1}\geq P^{min}_k$ and $P^{max}_{k+1}\geq P^{max}_k$, then
\beqs 
P^{min}_{k+2}&=&c^TP^{max}_{k+1} \; \geq \;  c^TP^{max}_{k} \; = \; P^{min}_{k+1}, 
\enqs  
and
\beqs 
P^{max}_{k+2}&=&c^TP^{min}_{k+2}+(1-c^T)\frac{\bc^0}{1-c} 
\; \geq \; 
c^TP^{min}_{k+1}+(1-c^T)\frac{\bc^0}{1-c} \; = \; P^{max}_{k+1},
\enqs
which concludes the induction. By \eqref{eq-framing-PMKV}, we deduce that $P^{min}$ and $P^{max}$ are two bounded increasing sequences, which thus converge to $P^{min}_\infty$ and $P^{max}_{\infty}$. The identity $P^{min}_k\leq P^{max}_{k}$, which holds true by definition, implies that $P^{min}_\infty\leq P^{max}_{\infty}$. Finally, by \eqref{Pkminmax}, we have
\beqs 
P^{min}_{k+1}
&=&c^{2T}P^{min}_k + c^T(1-c^T)\frac{\bc^0}{1-c},
\enqs 
and thus in the limit $k\rightarrow \infty$, we obtain
\beqs 
P^{min}_{\infty}&=& c^{2T}P^{min}_{\infty} + c^T(1-c^T)\frac{\bc^0}{1-c}
\enqs 
which yields
\beqs 
P^{min}_{\infty} &=& \frac{c^T}{1+c^T}\frac{\bc^0}{1-c}.
\enqs 
Likewise, we obtain
\beqs 
P^{max}_{\infty} &=&  \frac{1}{1+c^T}\frac{\bc^0}{1-c}.
\enqs

\printbibliography

\end{document}